\documentclass[12pt,reqno]{amsart}
\usepackage{amssymb}

\newcommand{\threecase}[7]{#1 \begin{cases} #2 & \text{{\rm #3}}\\ #4
&\text{{\rm #5}}\\ #6 & \text{{\rm #7}} \end{cases}   }

\newcommand{\sym}{\text{\rm sym}}

\newcommand{\supp}{{\rm supp}}

\newcommand{\mattwo}[4]
{\left(\begin{array}{cc}
                        #1  & #2   \\
                        #3 &  #4
                          \end{array}\right) }

\newcommand\be{\begin{equation}}
\newcommand\ee{\end{equation}}
\newcommand\bp{\begin{proof}}
\newcommand\ep{\end{proof}}

\newcommand\bea{\begin{eqnarray}}
\newcommand\eea{\end{eqnarray}}
\newcommand\bml{\begin{multline}}
\newcommand\eml{\end{multline}}
\newcommand\bal{\begin{align}}
\newcommand\eal{\end{align}}
\newcommand\bi{\begin{itemize}}
\newcommand\ei{\end{itemize}}
\newcommand\ben{\begin{enumerate}}
\newcommand\een{\end{enumerate}}
\newcommand\bc{\begin{center}}
\newcommand\ec{\end{center}}
\newcommand\ba{\begin{array}}
\newcommand\ea{\end{array}}

\newtheorem{thm}{Theorem}[section]

\newtheorem{cor}[thm]{Corollary}
\newtheorem{lem}[thm]{Lemma}

\newtheorem{exa}[thm]{Example}
\newtheorem{defi}[thm]{Definition}

\newtheorem{rek}[thm]{Remark}

\newcommand{\tf}{\widetilde{f}}
\newcommand{\tg}{\widetilde{g}}
\newcommand{\whf}{\widehat{f}}
\newcommand{\whg}{\widehat{g}}

\def\notdiv{\ \mathbin{\mkern-8mu|\!\!\!\smallsetminus}}
\newcommand{\R}{\ensuremath{\mathbb{R}}}
\newcommand{\A}{\ensuremath{\mathbb{A}}}
\newcommand{\C}{\ensuremath{\mathbb{C}}}

\newcommand{\Q}{\mathbb{Q}}

\newcommand{\ga}{\alpha}     
\newcommand{\gep}{\epsilon}   

\newcommand{\hphi}{\widehat{\phi}}  

\newcommand{\fof}{\frac{1}{4}}  
\newcommand{\foh}{\frac{1}{2}}  

\newcommand{\F}{\mathcal{F}} 
\newcommand{\G}{\mathcal{G}}

\newcommand{\hfo}{\widehat{f_1}}
\newcommand{\hft}{\widehat{f_2}}

\newcommand{\kkot}[1]{ \frac{\sin \pi {#1} }{\pi {#1} } }

\newcommand{\fn}{\mathcal{F}_N}
\newcommand{\gm}{\mathcal{G}_M}
\newcommand{\f}{\mathcal{F}}
\newcommand{\g}{\mathcal{G}}
\newcommand{\gln}{\text{{\rm GL}}_n}
\newcommand{\glt}{\text{{\rm GL}}_2}

\newcommand{\glm}{\text{{\rm GL}}_m}
\newcommand{\un}{\text{{\rm U}}}
\newcommand{\sy}{\text{{\rm Sp}}}
\newcommand{\soe}{\text{{\rm SO(even)}}}
\newcommand{\soo}{\text{{\rm SO(odd)}}}
\newcommand{\so}{\text{{\rm O}}}
\newcommand{\ntg}{\text{{\rm NT-good}}}

\newcommand{\afn}{|\fn|}

\newcommand{\tafn}{\frac2{\afn}}
\newcommand{\agm}{|\gm|}

\newcommand{\bfA}{\mathbb{A}}
\newcommand{\bfQ}{\mathbb{Q}}
\newcommand{\bfR}{\mathbb{R}}
\newcommand{\bfC}{\mathbb{C}}
\newcommand{\bfZ}{\mathbb{Z}}
\newcommand{\cH}{\mathcal{H}}
\newcommand{\WR}{\mathcal{W}_\bfR}
\newcommand{\GR}{\Gamma_\bfR}
\newcommand{\GC}{\Gamma_\bfC}

\DeclareMathOperator{\sgn}{sgn}
\DeclareMathOperator{\GL}{GL}

\numberwithin{equation}{section}

\begin{document}

\title[The effect of convolving on the underlying group symmetries]
{The effect of convolving families of $L$-functions\\ on the underlying group
symmetries}

\author{Eduardo Due\~nez} \email{eduenez@math.utsa.edu}
\address{Department of Mathematics, The University of Texas at San Antonio,
San Antonio, TX   78249}

\author{Steven J. Miller}
\email{Steven.J.Miller@williams.edu} \address{Department of Mathematics and Statistics, Williams College, Williamstown, MA 01267}


\thanks{We thank Jim Cogdell, Sol Friedberg, Gergely Harcos, Jeff
Hoffstein, Wenzhi Luo, Stephen D. Miller, Steve Rallis, Ze\'{e}v
Rudnick, Peter Sarnak, Eitan Sayag and Joe Silverman for many
enlightening conversations, and Owen Barrett for catching a typo in an earlier draft. The first named author was partly
supported by EPSRC grant N09176; the second named author was
partially supported by NSF grant DMS-0600848.}

\date{\today}

\begin{abstract}
  Let $\{\fn\}$ and $\{\gm\}$ be families of primitive automorphic
  $L$-functions for $\gln(\bfA_\bfQ)$ and $\glm(\bfA_\bfQ)$,
  respectively, such that, as $N,M \to \infty$, the statistical
  behavior ($1$-level density) of the low-lying zeros of $L$-functions
  in $\fn$ (resp., $\gm$) agrees with that of the eigenvalues near $1$
  of matrices in $G_1$ (resp., $G_2$) as the size of the matrices tend
  to infinity, where each $G_i$ is one of the classical compact groups
  (unitary $\un$, symplectic $\sy$, or orthogonal~$\so, \soe,
  \soo$). Assuming that the convolved families of $L$-functions $\fn
  \times \gm$ are automorphic, we study their $1$-level density. (We
  also study convolved families of the form $f\times \gm$ for a
  fixed~$f$.)  Under natural assumptions on the families (which hold
  in many cases) we can associate to each family $\mathcal{L}$ of
  $L$-functions a symmetry constant $c_\mathcal{L}$ equal to $0$
  (resp., $1$ or $-1$) if the corresponding low-lying zero statistics
  agree with those of the unitary (resp., symplectic or orthogonal)
  group. Our main result is that $c_{\f\times\g} = c_\f \cdot c_\g$:
  the symmetry type of the convolved family is the product of the
  symmetry types of the two families. A similar statement holds for
  the convolved families $f\times\gm$. We provide examples built from
  Dirichlet $L$-functions and holomorphic modular forms and their
  symmetric powers. An interesting special case is to convolve two
  families of elliptic curves with positive rank. In this case the symmetry
  group of the convolution is independent of the ranks, in accordance
  with the general principle of multiplicativity of the symmetry
  constants (but the ranks persist, before taking the limit
  $N,M\to\infty$, as lower-order terms).
\end{abstract}

\maketitle


\section{Introduction}

\subsection{Preliminaries }\ \\

Many questions in number theory, such as the study of the density of
the primes or properties of class numbers, can be related to
understanding the distribution of zeros of $L$-functions. In the early
1970's Dyson and Montgomery~\cite{Mon} discovered the agreement
between the pair correlation of zeros of the Riemann zeta-function
$\zeta(s)$ and of eigenvalues of matrices in the Gaussian Unitary
Ensemble~(GUE). Two decades later Katz and Sarnak~\cite{KaSa1,KaSa2}
offered deeper insight into the connection between zero and eigenvalue
statistics by studying \emph{families} of $L$-functions.  Ever since,
random matrix theory \cite{CFKRS,KaSa2,KeSn,ILS} has enjoyed
remarkable success at modeling and predicting the behavior of
$L$-functions.

Various pairs of eigenvalue/zero statistics can be shown, or at least
are conjectured, to be in perfect agreement. Among the early such
statistics studied were $n$-level correlations and
nearest-neighbor spacings~\cite{Hej,Mon,Od1,Od2,RS}. These statistics
pertain to the whole (infinite) sequence of critical zeros of a
\emph{single} $L$-function, and are shown to agree with the
corresponding statistic of $N\times N$ GUE matrices in the limit as
$N\to\infty$.  Neither of these statistics reveals anything about the
behavior of low-lying critical zeros of $L$-functions; that is, of
zeros near the arithmetically-crucial central point. The reason is
that those statistics are defined by averaging quantities defined
using a large (but finite) subset of the zeros, most of which will lie
high up on the critical line ---and thus the behavior of those few
zeros that lie near the central point is irrelevant in the limit as
the number of zeros used to compute the statistic tends to infinity.

The Katz-Sarnak philosophy has shifted the emphasis
to the study of families of $L$-functions and their low-lying
zeros, whose statistics (upon averaging over the family) are well
modeled by the statistics of eigenvalues close to~$1$ of random
matrices from the classical compact groups. In the function field case
these classical group statistics are explained by the monodromy group
of the family. For families of automorphic $L$-functions of number
fields the connection is quite a bit more mysterious. Usually, the
corresponding classical compact group is identified only by explicitly
computing zero statistics. Our goal in this paper is to allow
predicting the group attached to a ``convolved'' family assuming only
knowledge of the groups describing the zero statistics of the two
families being convolved. The relation turns out to be very simple to
describe and it will hopefully shed some light into the properties of
the (conjectural) correspondence between families of (number-field)
automorphic $L$-functions and classical groups.

We first describe the main statistic studied in this paper. In order
to break away from the universal global GUE statistics of the zeros of
a single $L$-function, and to understand the neighborhood of the
central point, we study the $n$-level density in a family of
$L$-functions; the latter is a local statistic involving only critical
zeros near the central point. Let $\mathcal{F} = \cup \mathcal{F}_N$
be a family of $L$-functions ordered by their conductors (for example,
$\mathcal{F}_N$ might be Dirichlet $L$-functions with conductor $N$ or
cuspidal newforms of weight $2$ and level~$N$) and write the zeros of
$L(s,f)$ as $1/2 + i \gamma_{j;f}$; assuming the General Riemann
Hypothesis (GRH), each $\gamma_{j;f} \in \R$. Given an $n$-variable
test function $\phi(t_1,\dots,t_n) =
\phi_1(t_1)\cdot\dots\cdot\phi_n(t_n)$ (where each $\phi_k$ is a
Schwartz function on $\R$), the \emph{$n$-level density} of
$\mathcal{F}$ is (by a slight abuse of language) the measure
on $\R^n$ with respect to which the integral of
$\phi$ is
\begin{eqnarray}\label{eq:defn1leveldensity}
D_{n,\mathcal{F}}(\phi) = \lim_{N\to\infty}\frac1{|\mathcal{F}_N|} \sum_{f\in\mathcal{F}_N} \sum_{j_1,\dots, j_n \atop j_i
\neq \pm j_k} \phi_1\left(\frac{\log
Q_f}{2\pi}\gamma_{j_1;f}\right)\cdots \phi_n\left(\frac{\log
Q_f}{2\pi}\gamma_{j_n;f}\right),
\end{eqnarray}
provided the limit exists, where $Q_f$ is the analytic conductor of
$L(s,f)$.\footnote{\label{fn:rescaling}It is a simple consequence of
  the Riemann-von Mangoldt zero-counting formula that the density of
  the zeros near the central point $s=1/2$ is roughly $(\log
  Q_f)/2\pi$, so the rescaled (``normalized'') imaginary parts
  $\gamma_{j_n;f}\cdot(\log Q_f)/2\pi$ have uniform (constant) density
  $1$ in the large-conductor limit. Thus, for fixed $A<B$, each
  $L(s,f)$ has roughly $B-A$ normalized zeros with imaginary parts
  in~$[A,B]$. Also, critical zeros not near $s=1/2$ (on a scale of
  $(\log Q_f)/2\pi$) are ``rescaled away to infinity'' in the
  large-conductor limit.}  For many families of $L$-functions
\cite{DM,FI,Gu,HR,HM,ILS,KaSa2,Mil2,Ro,Rub,Yo2} (and, conjecturally at
least, for any natural such family, in accordance with the Katz-Sarnak
philosophy), the $n$-level density coincides with that of the
normalized eigenvalues near $1$ of matrices in one of the infinite
families of classical compact Lie groups, in the limit as the size $N$
of the matrix goes to infinity. In the context of matrices from
classical Lie groups, the averaging over $\mathcal{F}_N$ in
equation~(\ref{eq:defn1leveldensity}) is replaced by averaging over
the whole group with respect to its natural (Haar) probability
measure ---hence the terminology of ``random matrices''.  The $n$-level
densities for different classical compact groups are distinct ---it is
this feature that allows ``breaking'' the universal GUE behavior
observed when considering global statistics such as $n$-level
correlations or neighbor spacings. This one-to-one correspondence
between (infinite families of) classical groups and their $n$-level
densities allows, at least conjecturally, to assign a definite
``symmetry type'' to each family of primitive $L$-functions. For
families of zeta or $L$-functions of curves or varieties over finite
fields, the corresponding classical compact group is determined by the
monodromy group of the family~\cite{KaSa1}. However, for families of
number-field automorphic $L$-functions there is no such thing as a
monodromy group and the underlying symmetry only manifests itself (in
our current understanding) through the zero statistics (although
function field analogues of number-field families often suggest what
the symmetry type should be). Our goal in this paper is to determine
the symmetry group of the convolution of two families of number-field
automorphic $L$-functions in terms of the symmetry groups of the
families being convolved together.

For families where the signs of the functional equations are all even
and there is not an obvious corresponding family with odd signs, a
folklore conjecture (see for example page 2877 of~\cite{KeSn}) stated
that the symmetry group should be symplectic. This was based on the
observation that $\text{{\rm SO(even/odd)}}$ symmetries in all known
examples arose from splitting orthogonal families according to the
sign of the functional equation, while symplectic symmetries arose
from a family with all even sign and no corresponding family with odd
signs. \emph{A priori} the symmetry type of a family with all
functional equations even is either symplectic or $\soe$. In \cite{DM}
we studied the family $\{L(s,\phi\times\sym^2f)\}$, where $\phi$ is a
fixed even Hecke-Maass eigenform on the modular group and $f$ ranges
over weight-$k$ full-level Hecke cusp forms; see \cite{LS} for
applications of this family. All $L(s,\phi\times\sym^2f)$ have even
sign, and this family does not arise from splitting sign within an
orthogonal family. In \cite{DM} it is shown (via $1$- and $2$-level
densities) that the symmetry type agrees only with $\soe$, thus
disproving the folklore conjecture mentioned above.

As a consequence of the counterexample to the folklore conjecture, the
theory of low-lying zeros is more than just a theory of the signs of
functional equations. By analyzing Rankin-Selberg convolutions of
$\glt$ $L$-functions (and some of their lifts), we are led to
attaching a symmetry constant $c_{\mathcal{F}}$ to each family
$\mathcal{F}$ of $L$-functions. This constant depends only on the
second moment (i.~e., the average over the family) of the Satake parameters at
each unramified prime. In all the cases investigated, the average
is~$0$ (resp., $1$ or~$-1$) if the family has unitary (resp.,
symplectic or orthogonal) symmetry. We are ready to set some notation
and describe our main result, namely that in many cases the symmetry
constant of the convolution of two families is the product of their
symmetry constants.

\bigskip
\subsection{$n$-Level Densities and $\ntg$ Families}\label{subsec:RMTgood}\ \\

We list four desirable properties for a family of primitive
$L$-functions to have; we call a family satisfying these properties $\ntg$. These properties are inspired by the families that have been
successfully investigated to date, and codify the conditions for
which we can calculate the $1$-level (and sometimes even the
$n$-level) densities for a family of $L$-functions. Though we could
replace some of the bounds with slightly weaker conditions, these are
the conditions that are met in practice.


\begin{defi}[$\ntg$]\label{defntg} Let $\phi$ be an even Schwartz
  test function such that ${\rm supp}(\hphi) \subset (-\sigma,
  \sigma)$ for some $\sigma > 0$. A family $\f$ of primitive
  automorphic $L$-functions for $\gln(\bfA_\bfQ)$ is $\ntg$ with
  symmetry constant $c_\f$ if $\f$ is a disjoint union of finite sets
  $\fn \subset \f$ such that, as $N \to \infty$:

\ben

\item \underline{\emph{Cardinality:}} \\
     {\rm  (i) Bounded multiplicities:} Members of a family can occur multiple times, say $f\in\f_N$ occurs $\mu_f$ times; however, we assume the multiplicities are bounded by a universal constant, independent of $N$: $\mu_f \le \mu_\f$.\\
     {\rm (ii) Size of the family:}  $\afn \to \infty$, where each member is counted with its multiplicity: $|\f_N| = \sum_{f\in\f_N} \mu_f$.\\

\item \underline{\emph{Conductors:}} The analytic log-conductors of
  $f\in\fn$ are essentially constant: say, $\log Q_f = \log R_N +
  o(\log R_N)$ for all $f\in\fn$ and some sequence $\{R_N\}$.  Further, there exists $\delta_0, \delta_0'>0$ such that
  $\afn^{\delta_0} \ll R_N \ll \afn^{\delta_0'}$. (In particular,
  $R_N\to\infty$ as $N\to\infty$.)\\

\item \underline{\emph{Sums over primes and squares of primes:}}%
  \footnote{The numbers $b_f(n)$ are the Dirichlet coefficients of
    the logarithmic derivative:\\ $L'(s,f)/L(s,f)=\sum_{n=1}^\infty
    b_f(n)/n^s$, cf., Definition~\ref{def:b} below.}\\
  {\rm (i)  Prime sums:} For some $r_{\mathcal{F}} \ge 0$,
  \be\label{eq:onelevelsumbfpb} -2 \sum_p \frac{1}{\sqrt{p}}
  \frac{\log p}{\log R_N}\ \hphi\left(\frac{\log p}{\log R_N}\right)
  \frac1{|\fn|} \sum_{f\in \fn} b_f(p)\ =\ r_\F \phi(0) +
  o\left(1\right); \ee we call $r_\F$ the \emph{rank of the
    family}. Often \eqref{eq:onelevelsumbfpb} is satisfied because
  $\exists \delta_1 > 0, \mu_1, r_{\mathcal{F}} \ge 0$ such that
  \bea\label{eq:onelevelsumbfp} \frac1{|\mathcal{F}_N|}\sum_{f\in\fn}
  b_f(p)\ = \ -\frac{r_{\mathcal{F}}}{\sqrt{p}} +
  O\left(\afn^{-\delta_1}p^{\mu_1}\right). \eea It also suffices for
  this to hold for almost all primes,
  provided the contribution from the bad primes is negligible.  \\

  \noindent {\rm (ii) Prime-square sums:} For some $c_\f \in \{-1,0,1\}$,
  \be\label{eq:3iitwo}
   -2 \sum_p \frac1p \frac{\log p}{\log R_N}
  \hphi\left(2\frac{\log p}{\log R_N}\right) \frac1{|\fn|} \sum_{f\in\fn}
  b_f(p^2) \ = \ -c_\f \frac{\phi(0)}2 + o(1).
  \ee Often \eqref{eq:3iitwo} is satisfied because $\exists \delta_2
  > 0, \mu_2 \ge 0, c_\f \in \{-1,0,1\}$ such that
  \bea\label{eq:3iione}
  \frac1{|\mathcal{F}_N|}\sum_{f\in\fn} b_f(p^2)\  =\ c_\f +
  O(\afn^{-\delta_2}p^{\mu_2}).
  \eea We call $c_\f$ the \emph{symmetry constant} of the family.\bigskip

\item \underline{\emph{Error terms:}} We have
  \be\label{eq:nuge3inftypsum} \frac1{|\fn|} \sum_{f\in\fn} \sum_{p}
  \sum_{\nu = 3}^\infty \frac{b_f(p^\nu)}{p^{\nu/2}} \frac{\log
    p}{\log R_N} \hphi\left(\nu\frac{\log p}{\log R_N}\right) \ = \
  o(1). \ee Many estimates on $b_f(p^{\nu})$ imply
  \eqref{eq:nuge3inftypsum}; we give three natural ones:\\
  {\rm (i) Ramanujan Conjecture:} For all $p$, $|\ga_{f,j}(p)| \le 1$
  implying $b_f(p^{\nu}) = O(1)$. \\
  {\rm (ii)} $\exists \mu_3(\nu) > 0$ with $\lim_{\nu\to\infty} \mu_3(\nu) =
  \mu_3 > 0$ such that \be b_f(p^\nu)\ \ll\
  \frac{p^{\nu/2}}{p^{1+\mu_3(\nu)\nu}}.\ee (iii) $\exists \delta_3 >
  0, \mu_3 \ge 0$ such that \be\sum_{f\in\fn} b_f(p^\nu)\ \ll\
  \afn^{1-\delta_3}p^{\mu_3}.\ee \een
\end{defi}

Condition~(1) ensures we have enough $L$-functions for
averaging, and when we convolve two families, it is a needed ingredient in controlling the contribution of imprimitive $L$-functions.\footnote{We could weaken our assumptions and allow $\mu_f \le |\f_N|^{1-\gep}$ for some $\gep > 0$. In the applications we have in mind, all multiplicities are bounded, so for ease of exposition we assume the multiplicities are bounded. We comment on this further in Corollary \ref{cor:cortomainthmonconvolvingfamilies}.} Condition~(2) allows us to handle the conductors,
and ensures that the number of $L$-functions is at least a power of
the analytic conductor; this is often needed in averaging to show
certain terms are small. Condition~(4) allows us to ignore
the contributions from $\nu\ge 3$ in the explicit formula
\eqref{eq:14}. The point of condition~(4.ii) is that eventually (for
$\nu$ large) we have $b(p^\nu)/p^{\nu/2} \ll
1/p^{1+\mu_3\nu}$, and this will be summable over $\nu$ and
$p$ (the small $\nu \ge 3$ terms can be handled individually by our
assumptions); an alternate bound where the cancelation
comes not from each individual $L(s,f)$ but rather from averaging
over the family is given by condition~(4.iii).

Condition~(3) is the interesting one, especially~(3.ii). It is
here that we see family-dependent behavior. In order to use the
Explicit Formula~\eqref{eq:14} successfully, we need to be able to determine family
averages of $b_f(p)$ and $b_f(p^2)$.  Condition~(3) holds in all the families of
$L$-functions studied to date
\cite{DM,FI,Gao,Gu,HM,HR,ILS,Mil1,Mil2,RR,Ro,Rub,Yo2};
further $r_\F$ is zero except for families of elliptic curves with positive
rank. The main term of the family averages of $b_f(p^2)$ do not depend
on $r_\F$, which surfaces only in the averages of $b_f(p)$.
Condition~(3) holds in the form~\eqref{eq:3iione} in all families studied to date
except for one-parameter families of elliptic curves with constant
$j$-invariant. (Michel~\cite{Mic} proved that~\eqref{eq:3iione} holds for
one-parameter families of elliptic curves with non-constant
$j$-invariant. If the $j$-invariant is constant one can often show by
direct computation that either \eqref{eq:3iitwo} holds or
\eqref{eq:3iione} holds on average; see~\cite{Mil1,Mil2}.)

We conclude this subsection by listing some $\ntg$ families with \emph{constant}
analytic conductors in each $\fn$.\\

\noindent Unitary \bi\item \{$L(s,\chi): \chi$ a non-trivial Dirichlet
character of prime conductor $m$\}, $m\to\infty$ (see~\cite{HR}).\\
\ei

\noindent Symplectic \bi \item $\{L(s,\chi_d): d$ ranges over subsets of fundamental discriminants in $[N,2N]\}$, $N\to\infty$ (see \cite{Gao,HR,Mil5,Rub}). \item \{$L(s,\sym^r f):$ $r$ even and $f$
ranges over weight-$k$ full-level cusp forms\}, $k\to\infty$ (see
\cite{Gu,ILS}). \item \{$L(s,\phi\times f): \phi$ a fixed Maass form
and $f$ ranges over weight-$k$ full-level cusp forms\}, $k\to\infty$
(see \cite{DM}). \item $\{L(s,\psi):\psi$ a character of the
ideal class group of the imaginary quadratic field
$\mathbb{Q}(\sqrt{-D})$ with $D>3$ square-free and congruent to 3
modulo 4$\}$ (see \cite{FI}). \\ \ei

\noindent Orthogonal \bi \item \{$L(s,f): f$ ranges over weight-$k$
level-$N$ cuspidal newforms with $k$, $N$ or both tending to
infinity\}; if we split by sign of the functional equations we get
$\soe$ or $\soo$ symmetry (\cite{ILS,Mil7,RR,Ro} for the $1$-level and
\cite{HM} for the $n$-level density). \item \{$L(s,\phi\times\sym^2
f): \phi$ a fixed Maass form and $f$ ranges over weight-$k$ full-level
cusp forms\}, $k\to\infty$ (see \cite{DM}). \item
\{$L(s,\sym^r f):$ $r$ odd and $f$ ranges over weight-$k$ full level
cusp forms\}, with $\so$ symmetry if $r\equiv 1,5 \bmod 8$, $\soe$
symmetry for $r \equiv 7\bmod 8$ and $\soo$ symmetry for $r\equiv 3
\bmod 8$, $k\to\infty$ (see \cite{Gu}).\\ \ei

With some more work, families with monotone increasing conductors
can be handled. This allows us to add an entry to
each list. For unitary families we may consider non-primitive
Dirichlet characters with square-free conductor (see \cite{Mil8}).
For symplectic we may consider primitive quadratic Dirichlet
characters (see \cite{Rub}). For orthogonal families we may consider
one-parameter (see \cite{Mil2}) or two-parameter (see \cite{Yo2})
families of elliptic curves. A generic one-parameter family should
have rank~$0$ and equidistribution of signs of functional equations,
giving $\so$ symmetry; however there are numerous families with
positive rank, as well as constant sign families (see \cite{Mil2}
for exact statements and details).

\bigskip
\subsection{Main Results}\ \\

\emph{We adopt the following convention throughout this paper: If $\f$
  and $\g$ are two families of unitary automorphic cuspidal
  representations of ${\rm GL}_n(\bfA_\bfQ)$ and ${\rm
    GL}_m(\bfA_\bfQ)$ with trivial central character, then by $\f
  \times \g$ we mean the set of all the (conjectural) Rankin-Selberg
  automorphic representations $f\times g$ of ${\rm
    GL}_{mn}(\bfA_\bfQ)$, where $f \in \f$, $g \in \g$. (Here the
  $f\times g$ are counted with multiplicity $\mu_f\mu_g$.)  For
  every purpose in this paper, this is equivalent to considering $\f$
  and $\g$ to be families of automorphic $L$-functions and
  $\f\times\g$ consists of the Rankin-Selberg convolution
  $L$-functions $L(s,f\times g)$.}

\bigskip

We occasionally remind the reader of this convention. We need some
control over the number of pairs $(f,g)$ where $g$ is the
contragredient of~$f$.  In this case the convolved $L$-function
$L(f\times g, s)$ is imprimitive, and thus its zeros are the
superposition of the zeros of at least two primitive $L$-functions. In
most cases the number and contribution of these imprimitive
$L$-functions to the $1$-level density is negligible.

\begin{defi}[Symmetry constant, family constant]
  We denote the \emph{symmetry constant} of the family $\f$ by
  $c_{\mathcal{F}}$. It equals $0$ (resp., $1$ or $-1$) if the
  $1$-level density of the family agrees with unitary (resp.,
  symplectic or any of the three orthogonal groups: $\so$, $\soe$ or
  $\soo$).  As the three orthogonal groups all have $c_{\mathcal{F}} =
  -1$, to distinguish them we set $\gep_\mathcal{F}$ equal to $0$
  (resp., $1$ or $-1$) if $\mathcal{F}$ has half of the signs of its
  functional equation even (resp., all signs even or odd); if
  $\mathcal{F}$ is not associated to an orthogonal group, then
  $c_\mathcal{F}$ alone determines the group and we simply put
  $\gep_\mathcal{F} = 0$. Finally, $r_{\mathcal{F}}$ denotes the rank
  of $\mathcal{F}$; except for families of elliptic curves, all other
  known families have $r_\mathcal{F}=0$. We call $\widetilde{c}_\F =
  (c_\F,\gep_\F,r_\F)$ the \emph{family constant} of~$\f$. \end{defi}

\begin{thm}\label{thm:ntgRSconv}
  Let $\f$ and $\g$ be $\ntg$ families of unitary automorphic cuspidal
  representations of ${\rm GL}_n(\bfA_\bfQ)$ and ${\rm
    GL}_m(\bfA_\bfQ)$ with trivial central character, with symmetry
  constants $c_\f$ and $c_\g$. Assume $\f \times \g$ is an $\ntg$
  family. Then the family $\f \times \g$ (which is the limit of $\fn
  \times \gm$, where $N$ and $M$ tend to infinity together) has
  symmetry constant \be c_{\f \times \g} \ = \ c_\f \cdot c_\g. \ee


If the family constants are
  $\widetilde{c}_\F = (c_\F,\gep_\F,r_\F)$ and $\widetilde{c}_\G =
  (c_\G,\gep_\G,r_\G)$ then the new family constant is
  $\widetilde{c}_{\F\times \G} = (c_\F \cdot c_\G,\gep_{\F\times
    \G},0)$.
\end{thm}

\begin{rek} Note that the ranks of the two families do not enter in the
  determination of the classical compact group associated to
  $\F\times\G$; the new family has rank~$0$ (in the sense of
  Definition~\ref{defntg}). Determining the distribution of signs of
  the functional equations of $\F\times\G$ is often an involved
  calculation depending on fine properties of the two families;
  however, it is not needed if we merely wish to classify the symmetry
  as unitary, symplectic or (non-specific) orthogonal. \end{rek}

\begin{rek} It is worthwhile to emphasize the meaning of the
  above theorem. Our notion of a symmetry constant arises from
  $1$-level density expansions, though we expect to see the same
  correspondence for any statistic ($n$-level correlations, central
  values, moments, $\dots$). Thus, an alternate way to phrase our
  results is that the $1$-level density of the convolution (as the
  conductors tend to infinity) agrees with the scaling limit of either
  unitary, symplectic or orthogonal matrices (depending on the value
  of the constant).
\end{rek}

\begin{rek} In the proof of Theorem \ref{thm:ntgRSconv}, our goal is to calculate the $n$-level densities for sufficiently large support to uniquely determine the corresponding symmetry group. A more involved argument could increase the support in some of our examples, but not far enough to see new features; we therefore content ourselves with giving the more general argument. We describe two examples in great detail which illustrate the technicalities that must be surmounted to prove the convolved families are $\ntg$; the first is symmetric powers of modular forms in \S\ref{sec:symmpowers}, and the second is families of elliptic curves in \S\ref{sec:ellcurves}. \end{rek}

\begin{cor}\label{cor:cortomainthmonconvolvingfamilies} The results of
  Theorem \ref{thm:ntgRSconv} still hold if we weaken the Bounded
  Multiplicities condition (1.i) in Definition~\ref{defntg}: Instead
  of assuming $\mu_f$ and $\mu_g$ are uniformly bounded, it suffices
  to assume that \be \#\{(f,g): f \in \f_N, g \in \g_M, f =
  \widetilde{g}\}\ =\ O\left(|\f_N|^{1-\delta} |\g_M| + |\f_N|
    |\g_M|^{1-\delta}\right) \ee for some $\delta > 0$ (in other
  words, that there is a power savings in the number of pairs where
  a~$g$ is the contragredient of an~$f$ ---these lead to imprimitive
  $L(s, f\times g)$).
\end{cor}

Instead of convolving $\f$ and $\g$, we can instead fix an
$f\in\f$ and consider the family $f\times\g$ obtained by taking the
limit as $M \to \infty$ of $f\times\gm$.

\begin{thm}\label{thm:twistfixedform}
  Assume $\g$ and $f\times\g$ are $\ntg$ and that $\g$ satisfies
  \eqref{eq:onelevelsumbfp} and \eqref{eq:3iione}. The symmetry type
  of $f\times \g$ is controlled by the following two pieces of input:
  $c_\g$ and $b_f(p^2)$. If $f$ is a Dirichlet character, holomorphic
  cusp form or Maass form then we may associate a symmetry constant
  $c_f$ to $f$ such that $c_{f\times\g} = c_f \cdot c_\g$. In
  particular, we have

\ben

\item if $f$ is a quadratic Dirichlet character then $f \times \G$ has
  the same symmetry as $\G$, and if $f$ is a non-quadratic Dirichlet
  character then $f\times\g$ has unitary symmetry;

\item if $\g$ has unitary (resp., symplectic, orthogonal) symmetry,
  then $f\times\g$ has unitary (resp., orthogonal, symplectic)
  symmetry if $f$ is a Hecke holomorphic or Maass form. \een
\end{thm}

\begin{rek} If instead \eqref{eq:onelevelsumbfpb} and
  \eqref{eq:3iitwo} hold then the result is probably still true (it
  can be shown in special cases by partial summation); in general,
  though, more detailed knowledge about sums of the coefficients of
  $L(s,f)$ will be needed. \end{rek}

\begin{rek}
The universality in Theorems \ref{thm:ntgRSconv}
and~\ref{thm:twistfixedform} is reminiscent of that found by Rudnick
and Sarnak~\cite{RS}, where the universality in the $n$-level
correlations of primitive automorphic cuspidal $L$-functions is
related to universality in the second moments of
the Fourier coefficients $a_\pi(p)$.\end{rek}

An especially interesting case is when at least one of the two
families is a one-parameter family of elliptic curves over $\Q(T)$
with positive rank $r_{\mathcal{F}}$. Miller~\cite{Mil2} showed that
the $1$- and $2$-level densities of zeros of these families agree with
those of \emph{subgroups} of the orthogonal group (in many cases
unconditionally, in other cases assuming standard conjectures; see
\cite{Yo2} for similar results involving special two-parameter
families). As the conductors tend to infinity, the random matrix
ensemble modeling this situation (as $N\to\infty$) is
\be\label{eq:scalingindepmodel}
\left\{\mattwo{I_{r_{\mathcal{F}}\times r_{\mathcal{F}}}}{}{}{g}, \ g
  \in \mathcal{C}\right\}, \ee where $I_{r_{\mathcal{F}}\times
  r_{\mathcal{F}}}$ is the $r_{\mathcal{F}}\times r_{\mathcal{F}}$
identity matrix and $\mathcal{C}$ is $\so(N)$ (resp., $\text{SO}(2N)$
or $\text{SO}(2N+1)$) if half the signs of the functional equation are
even (resp., all or none); the correct model is not known for finite
conductors (but see \cite{Mil4} for numerical investigations of zeros
near the central point). Indeed, by Silverman's specialization theorem
and the Birch and Swinnerton-Dyer conjecture, for all $t$ sufficiently
large each elliptic curve has at least $r_{\mathcal{F}}$ zeros at the
central point; moreover, the ensemble~(\ref{eq:scalingindepmodel})
models these zeros as independent from the remaining others. This
independence is in agreement with function-field analogues. We shall
see in Theorem~\ref{thm:maintwistecf} that if one convolves two
families of elliptic curves with positive rank then, to first order,
one does not see any effects of this rank in the symmetry group of the
new family! What this means is that the rank parameter $r$ of the new
family is zero as defined by Condition~(3) of Definition~1.1. The
ranks of the convolved families appear \emph{only} as a lower-order
correction term that is unfortunately difficult to isolate since it is
smaller than the bounds we can prove for the error terms (though,
conjecturally, it is larger than the actual bounds for these terms and
should, in principle, be detectable).  In this regard our results are
similar to Goldfeld's~\cite{Go}; he considered twists of a fixed
elliptic curve by quadratic Dirichlet characters and conjectured that
the new family's rank is independent of the rank of the fixed
elliptic curve.

In \S\ref{sec:central-1-2-level} and \S\ref{sec:random-matrix-theory}
we review the needed results from number theory and random matrix
theory. We discuss the properties and consequences of being an $\ntg$
family of $L$-functions in \S\ref{sec:RMTgood} and then in
\S\ref{sec:proofmainresults} prove Theorems \ref{thm:ntgRSconv} and
\ref{thm:twistfixedform}. We then give examples of families where
these conditions are met: Convolving families of holomorphic cusp
forms in Example \ref{exa:holocusp}, symmetric powers of holomorphic
cusp forms in \S\ref{sec:symmpowers}, and one-parameter families of
elliptic curves in \S\ref{sec:ellcurves}; these examples are all
independent of each other and may be read in any order.


\bigskip
\section{Number Theory Review}\label{sec:central-1-2-level}

We quickly review the notion of automorphic $L$-functions.  These are
the $L$-functions attached to automorphic representations
of~$\gln(\bfA_\bfQ)$. Our examples are built out of objects in ${\rm
  GL}_1$ (Dirichlet $L$-functions) and ${\rm GL}_2$ (Maass forms and
holomorphic modular forms, including those attached to elliptic
curves). We build more complicated $L$-functions by taking
Rankin-Selberg convolutions and other natural functorial operations
(e.~g., forming the symmetric square $L$-functions). These
constructions take us beyond ${\rm GL}_2$. It is impossible to cover
here but the barest facts about automorphic $L$-functions; see
\cite{Bor,Jac,JPS,RS} for more details.  We will focus on primitive
$L$-functions; these are attached to cuspidal representations and
cannot be further factored as products of other $L$-functions, hence
their critical zeros form an irreducible set in this sense.

Let $\pi = \hat\otimes_v \pi_v$ be a unitary irreducible cuspidal
automorphic representation of $\gln(\bfA_\bfQ)$ with trivial central
character. Here $v$ is either a prime $p$ or $\infty$, and each $\pi_v$ is
an irreducible admissible representation of $\bfQ_v$ (where
$\bfQ_\infty:=\bfR$). The finite part of the $L$-function attached
to $\pi$ is an Euler product \be L(s,\pi)\ =\ \prod_p L(s,\pi_p).
\ee Outside a finite set of primes, $\pi_p$ is unramified and \be
L(s,\pi_p) \ = \ \det(I - p^{-s} A_\pi(p))^{-1} \ = \ \prod_{j=1}^n
(1 - \ga_{\pi,j}(p)p^{-s}), \ee where $\{A_\pi(p)\} \in \gln(\C)$ is
a semi-simple conjugacy class parametrized by the eigenvalues
$\ga_{\pi,j}(p)$.  The Satake correspondence is the bijection
$A_\pi(p)\leftrightarrow\pi_p$ between semi-simple conjugacy classes
in $\gln(\bfC)$ and unramified irreducible admissible
representations of~$\gln(\bfQ_p)$.

The complex numbers $\{\ga_{\pi,j}(p)\}_{j=1}^n$ are called the Satake
parameters of $\pi_p$.  In the context at hand, the generalized
Ramanujan conjecture is the statement that $|\ga_{\pi,j}|=1$ at the
unramified places; at a ramified prime $p$ some of the
$\ga_{\pi,j}(p)$ may vanish.

\begin{defi}\label{def:b}
  For $\pi$ an automorphic representation, $p$ a prime and $\pi_p$
  with Satake parameters $\ga_{\pi,1}(p), \dots, \ga_{\pi,n}(p)$, we
  define, for $\nu=1,2,3,\dots$, \be b_\pi(p^\nu) \ :=\
  \ga_{\pi,1}(p)^\nu + \dots + \ga_{\pi,n}(p)^\nu. \ee
\end{defi}
With this definition one has $b_\pi(p^\nu)={\rm Trace}(A_\pi(p)^\nu)$
for unramified~$p$, and \be\frac{L'(s,\pi)}{L(s,\pi)} =
\sum_p\sum_{\nu=1}^\infty\frac{b_\pi(p^\nu)}{p^{\nu{s}}}.\ee

The archimedean $L$-factor associated to $\pi_\infty$ is of the form
\be L(s,\pi_\infty)\  =\  \prod_{j=1}^n \Gamma_\R(s+\mu_{\pi,j}),
\qquad \text{where $\Gamma_\R(s) = \pi^{-s/2}
\Gamma\left(\frac{s}2\right)$.} \ee The numbers
$\{\mu_{\pi,j}\}_{j=1}^n$ are analogs of the Satake parameters, and
the analog of the Ramanujan conjecture is in this case Selberg's
(generalized) eigenvalue conjecture, namely the statement that the
$\mu_{\pi,j}$ have non-negative real part.

We define the completed $L$-function by
\begin{equation}
\Lambda(s,\pi)\ :=\  N_\pi^{s/2} L(s,\pi_\infty) L(s,\pi),
\end{equation}
where $N_\pi$ is a positive integer called the arithmetic conductor.
We have the functional equation \be \Lambda(s,\pi)\ =\ \epsilon(\pi)
\Lambda(1-s,\widetilde{\pi}), \ee where $\widetilde{\pi}$ is the
contragredient of $\pi$ and $\epsilon(\pi)$ is a complex constant
such that $|\epsilon(\pi)|=1$.  In the self-dual case, when
$\pi\simeq\widetilde\pi$, the functional equation relates
$L(\cdot,\pi)$ to itself, and $\epsilon(s,\pi)$ equals $\pm1$.

For our applications, it is the analytic conductor (not the
arithmetic conductor) that is important for understanding the
behavior of the zeros near the central point. The two are related,
and the analytic conductor may be taken as \be\label{eq:defanalcond}
Q_\pi \ = \ \mu_{\pi,1}\ \cdots\ \mu_{\pi,n}\ N_\pi. \ee We use the
analytic conductor to rescale the low lying zeros, and then apply
the explicit formula to convert sums of an even Schwartz test
function over the zeros of the $L$-function to sums of the Fourier
transform of the test function evaluated at prime powers. For such
calculations, it is the logarithm of the analytic conductor that
normalizes the zeros; see for example section 4 of \cite{ILS}. In
some other papers our factors of $\mu_{\pi,j}$ are replaced with
$\mu_{\pi,j}'/2$. As we shall always be interested in situations
where the analytic conductors tend to infinity, both normalizations
lead to the same results. Note that we have $N_\pi^{s/2}$ in our
functional equation
---other authors sometimes write this factor as $(N_\pi')^s$,
which would lead to a factor of $(N_\pi')^2$ in the analytic
conductor.\\

Throughout the paper we make the following two assumptions, unless
specified.

\begin{itemize}
\item \emph{We assume the Generalized Riemann Hypothesis for all
    automorphic $L$-functions.} Thus we may write the non-trivial
  zeros as $\foh + i\gamma$ with $\gamma\in\R$, and the correct
  scaling for zeros near the central point is $\gamma \mapsto
  \tilde\gamma=\gamma \frac{\log Q_\pi}{2\pi}$ (low-lying
  $\tilde\gamma$'s have natural uniform density~$1$ ---see
  footnote~\ref{fn:rescaling}).  However, the results on $1$-level
  densities may be interpreted, and remain true, even when the
  $\gamma$ are allowed to be complex; see for example \cite{ILS}.  In
  other instances (such as in~\S\ref{eq:goodeqbfp2}), GRH is used to
  bound error terms and
  thus enters in the argument in a more essential manner.\\
\item \emph{We assume the Langlands functoriality conjectures for the
    automorphic representations under consideration.}  This is
  necessary in order to ensure that their attached automorphic
  $L$-functions have good analytic properties.  In some cases the
  analytic properties of an $L$-function are known even without
  knowledge of its automorphicity (e.g., for some symmetric-power
  $L$-functions attached to holomorphic modular
  forms~\cite{KiSh1,KiSh2,K}).  On the other hand, the automorphicity
  of all symmetric powers of an automorphic representation $f$ implies
  the Ramanujan-Selberg conjectures for $L(s,f)$, though bounds
  towards this goal are in some cases available
  unconditionally~\cite{K} (e.g.,
  Deligne's proof of Ramanujan for holomorphic modular forms).\\
\end{itemize}
\begin{rek}
  Automorphic $L$-functions associated to \emph{cuspidal}
  representations are \emph{primitive} in the sense that one cannot
  write $\Lambda(s,\pi)$ as
  $\Lambda(s,\pi_1)\Lambda(s,\pi_2)$. However, a general non-cuspidal
  automorphic $L(s,\pi)$ factors as a product of primitive ones and
  its critical zeros are clearly a union of the zeros of its primitive
  factors.  While the low lying zeros of each primitive factor will
  reveal a specific underlying symmetry (at least conjecturally), the
  low lying zeros of the imprimitive function will in general not
  correspond to a definite symmetry.  Here we are using the word
  ``symmetry'' in the sense of \S\ref{sec:random-matrix-theory}.
  However, even the assumption of functoriality does not ensure that
  lifts of cuspidal forms are cuspidal.  The simplest counterexample
  is the imprimitive $L$-function $L(s,\pi\times\widetilde\pi)$ where
  $\pi$ is a cuspidal automorphic representation of ${\rm
    GL}_n(\bfA_\bfQ)$, $n>1$, and $\widetilde\pi$ is its
  contragredient.%
  \footnote{If $\pi_1,\pi_2$ are automorphic unitary cuspidal
    representations as above, but not necessarily normalized to have
    trivial central character, then $L(s,\pi_1\times\pi_2)$ is
    imprimitive when
    $\pi_2\simeq\widetilde\pi_1\otimes|\det(\cdot)|^s$ is a twist of
    the contragredient of $\pi_1$.}  We will need to deal with this
  possibility on occasion.
\end{rek}

While we consider quite general families of $L$-functions, the
building blocks for examples which we can prove satisfy the necessary
conditions are (the automorphic representations attached to) either
Dirichlet characters or modular forms. For their corresponding
$L$-functions, classical summation formulas for Fourier coefficients
are available that make our approach tractable.


Let $\phi$ be an even Schwartz test function on~$\R$ whose Fourier
transform
\begin{equation}
  \label{eq:9}
  \widehat{\phi}(y)\ =\ \int_{-\infty}^\infty \phi(x) e^{-2\pi ixy}dx
\end{equation}
has compact support. Let $\mathcal{F}$ be a finite family of
$L$-functions satisfying GRH. For example, $\mathcal{F}$ might be the
set of all $L(s,\chi)$ with $\chi$ a non-trivial Dirichlet character
of conductor~$m$, and we would then investigate the limit as
$m\to\infty$. Other examples include weight-$k$ level-$N$ cuspidal
newforms (and let either $k$, or $N$, or both tend to infinity), as
well as one-parameter families of elliptic curves (where the parameter
$t$ varies over an interval $[N,2N]$, and then we let $N\to\infty$).


%

Consider a family $\mathcal{F}$ of $L$-functions $L(s,f)$ and denote
by $Q_f$ the analytic conductor of $L(s,f)$. Let $\mathcal{F}_N$ be
the finite subfamily of $\mathcal{F}$ consisting of those functions
with $Q_f = N$. Thus, $\mathcal{F} = \cup_N \mathcal{F}_N$. To study
the zeros of the functions in the family $\mathcal{F}$, we use the
Explicit Formula to convert sums over zeros to sums over primes. For
any $L$-function $L(s,f)$~\cite{ILS,RS}:
\begin{equation}
  \label{eq:14}
  \sum_\ell \phi\left(\gamma_{j;f}\frac{\log R}{2\pi}\right) \ = \
  \frac{A_f}{\log R}\ \hphi(0) -
  2\sum_{p}\sum_{\nu=1}^\infty\hphi\left(\frac{\nu\log p}{\log R}\right)
  \frac{b_f(p^\nu)\log p}{p^{\nu/2}\log R},
\end{equation}
where $A_f$ is an integral of gamma factors coming from the
functional equation of $L(s,f)$. We have
\begin{equation}
\label{eq:afqfrN} A_f\ =\
\log Q_f + o(\log Q_f),
\end{equation}
and the little-oh implicit constant often depends only on
$\mathcal{F}$ and not the individual $f$. We shall also consider
variations of the above family; for example, we may let
$\mathcal{F}_N$ be the set of $f$ in $\mathcal{F}$ with $N \le
Q_f \le 2N$. While the subject is considerably simplified if the
conductors in $\mathcal{F}_N$ are constant, monotonically increasing
conductors can be handled with additional work (see \cite{Mil2} for
details for families of elliptic curves).

After averaging over the family, the resulting sums are often
evaluated using the following consequence of the Prime Number Theorem:

\begin{thm}\label{thmprimesums} Let $\widehat{F}$ be an even Schwartz
  function of compact support. Then for any positive integer $\nu$,
  \be\label{eq:19} \sum_p \widehat{F}\left( \nu \frac{\log p}{\log R}
  \right) \frac{\log p}{\log R} \frac{1}{p} \ = \ \frac{1}{2\nu} F(0)
  + O\left( \frac{1}{\log R} \right). 
\ee
\end{thm}


\bigskip
\section{Random Matrix Theory Review}
\label{sec:random-matrix-theory}

Katz and Sarnak conjecture that to any infinite family $\mathcal{F}$
of $L$-functions one can associate one of the infinite families of
classical compact matrix groups (unitary, orthogonal, or symplectic),
say $G(\mathcal{F})$, such that the large-conductor statistics of
zeros near the central point for $L(s,f)$, $f\in\mathcal{F}$, agree
with those of the eigenvalues near $1$ of matrices in $G(\mathcal{F})$
(as the matrix size $N\to\infty$). Specifically, the $n$-level density
of rescaled critical zeros for the family $\mathcal{F}$ is the
function $W_{n,\mathcal{F}}$ (more precisely, the important
object is the measure $W_{n,\mathcal{F}}\,dx_1\dots dx_n$ on $\R^n$)
defined by its action on any test function
$\phi(x_1,\dots,x_n)=\phi_1(x_1)\cdots\phi(x_n)$ (where
$\phi_1,\dots,\phi_n$ are Schwartz functions) by: \bea
\label{eq:n-level-density}
D_{n,\mathcal{F}}(\phi)& \ =\ & \lim_{N\to\infty}
\frac{1}{|\mathcal{F}_N|}
\sum_{f\in \mathcal{F}_N} \sum_{\substack{j_1,\dots, j_n \\
    j_i \neq \pm j_k}} \phi_1\left(\gamma_{j_1;f}\frac{\log
    Q_f}{2\pi}\right)\cdots \phi_n\left(\gamma_{j_n;f}\frac{\log
    Q_f}{2\pi}\right) \nonumber\\ & = & \int \cdots \int
\phi_1(x_1)\cdots \phi_n(x_n)
W_{n,\mathcal{F}}(x_1,\dots,x_n)dx_1\cdots dx_n \nonumber\\ & = & \int
\cdots \int \widehat{\ \phi_1}(u_1)\cdots \widehat{\phi_n}(u_n)
\widehat{W}_{n,\mathcal{F}}(u_1,\dots,u_n)du_1\cdots du_n.  \eea (The
Fourier transform $\widehat W_{n,\mathcal{F}}$ is most often used in
proofs for technical reasons, e.~g., to use the Explicit
Formula~(\ref{eq:14}).) The Katz-Sarnak philosophy posits that the
$n$-level densities $W_{n,\mathcal{F}}$ agree with the $n$-level
densities $W_{n,G(\mathcal{F})}$ of the matrix group $G(\mathcal{F})$
associated to the family~$\mathcal{F}$. This philosophical
correspondence has been proved for many families when the Schwartz
test functions $\phi_i$ have Fourier transforms supported in a
sufficiently small neighborhood of zero. The $n$-level densities for
the classical compact groups are (see~\cite{KaSa1}):
\be\label{eqdensitykernels}\begin{array}{lcl}
  W_{n,\soe}(x) & \ = \ & \det (K_1(x_i,x_j))_{i,j\leq n} \\
  W_{n,\soo}(x) & = & \det (K_{-1}(x_i,x_j))_{i,j\leq n} +
  \sum_{k=1}^n
  \delta(x_k) \det(K_{-1}(x_i,x_j))_{i,j\neq k}\\
  W_{n,\so}(x) & = & \foh W_{n,\soe}(x) + \foh W_{n,\soo}(x) \\
  W_{n,\sy}(x) &=& \det (K_{-1}(x_i,x_j))_{i,j\leq n} \\
  W_{n,\un}(x) & = & \det (K_0(x_i,x_j))_{i,j\leq n},\end{array} \ee
where $K(y) = \kkot{y}$, $K_\epsilon(x,y) = K(x-y) + \epsilon K(x+y)$
for $\epsilon = 0, \pm 1$, and $\delta(u)$ is the Dirac Delta
functional.\footnote{While these determinant formulas hold for
  arbitrary support, in practice the resulting formulas for $n \ge 3$
  require some combinatorics when the support is large before
  agreement is seen with number theory. Hughes and Miller \cite{HM}
  derive an alternate formula for $n$-level statistics; while their
  formula holds for smaller support, in the range where it is
  applicable it facilitates comparisons with number theory.} The
Fourier transforms of the $1$-level densities are \be
\begin{array}{lcl}
\widehat{W}_{1,\soe}(u) & = & \delta(u) + \foh \eta(u) \\
\widehat{W}_{1,\soo}(u) & = & \delta(u) - \foh \eta(u) + 1
\\ \widehat{W}_{1,\so}(u) & = & \delta(u) + \foh
\\ \widehat{W}_{1,\sy}(u) & = & \delta(u) - \foh
\eta(u) \\ \widehat{W}_{1,\un}(u) & = & \delta(u),
\end{array}\ \ee where
\be \threecase{\eta(u) \ = \ }{1}{if $|u| < 1$}{\foh}{if $|u| =
  1$}{0}{if $|u| > 1$.} \ee When working with test functions $\phi$
whose Fourier transform $\widehat\phi$ is supported in a small
neighborhood of~$0$, it is still possible to distinguish between unitary,
symplectic and orthogonal $n$-level densities; however, as long as
$\widehat\phi$ is supported in $(-1,1)$, all three flavors of
orthogonal symmetry (even, odd, or full) agree: \be\label{eq:1ldwithtestfn}
\begin{array}{lcl}
\int \hphi(u) \widehat{W}_{1,\soe}(u)du & = & \hphi(u) + \foh \phi(0)\\
\int \hphi(u)\widehat{W}_{1,\soo}(u)du & = & \hphi(u) + \foh \phi(0)
\\ \int \hphi(u)\widehat{W}_{1,\so}(u)du & = & \hphi(u) +
\foh\phi(0) \\ \int \hphi(u)\widehat{W}_{1,\sy}(u)du & = & \hphi(u)
- \foh \phi(0) \\ \int \hphi(u)\widehat{W}_{1,\un}(u)du & = &
\hphi(u).
 \end{array}\ee Let ${\rm sign}(G) = 0$ (resp., $\foh$, $1$)
 for $G = \soe$ (resp., $\so$, $\soo$). For even functions
 $\phi(x_1,x_2)=\phi_1(x_1)\phi_2(x_2)$ such that
 $\hphi(u_1,u_2)=\hphi_1(u_1)\hphi_2(u_2)$ is supported in $|u_1| + |u_2| < 1$,
\begin{eqnarray}\label{eq:soetwo} \begin{array}{lcl}
\int \int \widehat{f_1}(u_1)\widehat{f_2}(u_2)
\widehat{W_{2,\mathcal{G}}}(u) du_1du_2 &= & \Big[\hfo(0) + \foh
f_1(0) \Big] \Big[\hft(0) + \foh f_2(0) \Big] \\ & & + \ 2 \int |u|
\hfo(u) \hft(u)du - 2 \widehat{f_1f_2}(0) - f_1(0)f_2(0)
\\ & & + \ {\rm sign}(\mathcal{G})f_1(0)f_2(0).
\end{array}
\end{eqnarray} Thus, for arbitrarily small support, the 2-level
density distinguishes the three orthogonal groups; see \cite{Mil1}
for the calculation.

In studying families of elliptic curves \cite{Mil2,Yo2}, often the
corresponding classical compact group is a subgroup of one of the
orthogonal groups. For one-parameter families of elliptic curves
over $\Q(T)$ with rank $r_{\mathcal{F}}$, as remarked in
\eqref{eq:scalingindepmodel}, the correct model as the conductors
tend to infinity appears to be \be
\left\{\mattwo{I_{r_{\mathcal{F}}\times r_{\mathcal{F}}}}{}{}{g}, \
g \in \mathcal{C}\right\}, \ee where $I_{r_{\mathcal{F}}\times
r_{\mathcal{F}}}$ is the $r_{\mathcal{F}}\times r_{\mathcal{F}}$
identity matrix and $\mathcal{C}$ is $\so$ (resp., $\soe$ or $\soo$)
if half the signs of the functional equation are even (resp., all or
none), though see \cite{Mil4} for a discussion of the behavior for
finite conductors. These $r_\f$ independent zeros replace $\hphi(u)
+ \foh \phi(0)$ with $\hphi(u) + \foh \phi(0) +
r_{\mathcal{F}}\phi(0)$ in the $1$-level density expansion, and
there is a similar modification in the $n$-level density.

Because of this effect of rank, we attach a family constant to each
family of $L$-functions $\mathcal{F}$: \be \widetilde{c}_{\mathcal{F}}
\ = \ (c_\mathcal{F},\gep_\mathcal{F},r_\mathcal{F}). \ee Here
$c_{\mathcal{F}}$ is the symmetry constant of the family, equal to $0$
(resp., $1$ or $-1$) if the family is unitary (resp., symplectic or
orthogonal); we call any subgroup of $\so$, $\soe$ or $\soo$
orthogonal. Since the three orthogonal groups all have
$c_{\mathcal{F}} = 1$, we set $\gep_\mathcal{F}$ equal to $0$ (resp.,
$1$ or $-1$) if $\mathcal{F}$ has half of the signs of its functional
equation even (resp., all signs even or odd); if $\mathcal{F}$ is not
associated to an orthogonal group, then $c_{\mathcal{F}}$ alone
determines the precise group and we define $\gep_\mathcal{F} =
0$. Finally, $r_{\mathcal{F}}$ denotes the rank of $\mathcal{F}$;
except for families of elliptic curves, all other known families have
$r_\mathcal{F}=0$.

Our main result (Theorem \ref{thm:ntgRSconv}) is that in order to
determine the symmetry of the the Rankin-Selberg convolution of two
$\ntg$ families, all that matters is $c_\mathcal{F}$ and
$c_{\mathcal{G}}$. Thus we may interpret the symmetry constant as a
convolution constant. Further, the new family has rank zero. This is
unfortunate, since otherwise this would allow constructing families of
$L$-functions with high central vanishing.


\bigskip
\section{$\ntg$ Families and $n$-Level Densities}\label{sec:RMTgood}

As a warm-up to proving our main theorems in
\S\ref{sec:proofmainresults}, in this section we investigate some
consequences of Definition~\ref{defntg} (NT-good).

It is worth commenting on the main terms in \eqref{eq:onelevelsumbfpb}
and~\eqref{eq:onelevelsumbfp}. Consider a one-parameter family of
elliptic curves over $\Q(T)$ with rank $r_\F$; $\mathcal{F}_N$ is
essentially just $\{E_t: t \in [N,2N]\}$.  Then $b_t(p) =
a_t(p)/\sqrt{p}$, where $a_t(p)$ are the coefficients of the
$L$-series of $L(s,E_t)$ (with functional equation $s \to 2-s$, so the
critical strip is $\Re s \in [0,2]$). Rosen and Silverman \cite{RoSi}
prove a conjecture of Nagao's (unconditionally if the elliptic surface
is rational; conditional on Tate's conjecture
otherwise): \begin{eqnarray} \lim_{X \to \infty} \frac{1}{X} \sum_{p
    \leq X} \frac{1}{p} \sum_{t=0}^{p-1} a_{t}(p) \log p\ =\ -r_\F.
\end{eqnarray} Thus the $a_t(p)$'s give the rank of
the family over $\Q(T)$. For many families of elliptic curves (see
\cite{ALM,Fe}), the main term of the average over $E_t\in\fn$ of
$a_t(p)/p$ is independent of $p$, and we have \bea
-\frac1{\afn}\sum_{E_t\in \mathcal{F}_N} \frac{b_t(p)}{\sqrt{p}} & \
= \ & -\frac1{\afn}\sum_{E_t\in \mathcal{F}_N}
\frac{a_t(p)}{p}\nonumber\\ &\ = \ & \frac1{\afn} \left[-
\frac{\afn}{p} \sum_{t\bmod p}
\frac{a_t(p)}{p} \right] + O\left(\frac{\sqrt{p}}{\afn}\right) \nonumber\\
& = & \frac{r_\F}{p} + O\left(\frac{\sqrt{p}}{\afn}\right).\eea If
we have such a family we use \eqref{eq:onelevelsumbfp}; if not, we
need to do a little more work and use \eqref{eq:onelevelsumbfpb} and
\eqref{eq:19}. The proofs follow similarly, the only real difference
being a partial summation on the primes to handle the test
functions.\\

\emph{For ease of exposition we concentrate on cases
where \eqref{eq:onelevelsumbfp} holds, and remark that similar
arguments handle the case when we have \eqref{eq:onelevelsumbfpb}.}


\begin{thm}\label{thm:1ntg} Let $\f$ be an $\ntg$ family of
  automorphic $L$-functions for $\gln$.  Then for even Schwartz test
  functions $\phi$ such that $\hphi$ is supported in a sufficiently
  small (but explicitly computable in terms of the constants
  $\delta_i, \mu_i$) neighborhood of~$0$, if $r_\f = 0$ then the
  1-level density of $\f$ agrees with unitary (resp., symplectic or
  orthogonal) if $c_{\mathcal{F}} = 0$ (resp., $1$ or $-1$); if $r_\F
  > 0$ then the corresponding classical compact group is modified by
  having an $r_\F \times r_\F$ identity matrix as
  in~\eqref{eq:scalingindepmodel}.
\end{thm}

\begin{proof} Using the explicit formula to calculate the 1-level
density, we have the expansion \bea D_{1,\fn}(\phi) &\ = \ &
\hphi(0)-\tafn \sum_{f\in\fn} \sum_{\nu=1}^\infty
\sum_{p=2}^{R_N^\sigma}  \frac{b_f(p^\nu) \log p}{p^{\nu/2}\log R_N}\ \hphi\left( \nu\frac{\log p}{\log R_N}\right) + o(1). \eea From
\eqref{eq:afqfrN}, the factor of $\hphi(0)$ above comes from the
constancy of the main term of the analytic conductors (and an
analysis of the $\Gamma$-factor terms; in fact, this is what we use
to determine $R_N$); the $o(1)$ term arises from the correction
factor in $\log Q_f = \log R_N + o(\log R_N)$. As our family is
$\ntg$, there is no contribution from $b_f(p^\nu)$ for $\nu \ge 3$
(either for all support, or for support sufficiently small). Thus
those terms may be absorbed into an error term.

We assume that \eqref{eq:onelevelsumbfp} holds; the case when
\eqref{eq:onelevelsumbfpb} holds follows similarly. The $\nu=1$
terms contribute \bea S_1 & \ = \ & -\tafn \sum_{f\in\fn}
\sum_{p=2}^{R_N^\sigma} \frac{b_f(p) \log p}{\sqrt{p}\log R_N}\ \hphi\left( \frac{\log p}{\log R_N}\right) \nonumber\\ & \ = \ &
2 \sum_{p=2}^{R_N^\sigma} \left[-\frac{1}{\afn} \sum_{f\in\fn}
b_f(p)\right] \frac{\log p}{\sqrt{p}\log R_N}\ \hphi\left(
\frac{\log p}{\log R_N}\right) \nonumber\\ & = & 2
\sum_{p=2}^{R_N^\sigma} \left[\frac{r_{\mathcal{F}}}{\sqrt{p}} +
O\left( \afn^{-\delta_1}p^{\mu_1}\right) \right] \frac{\log
p}{\sqrt{p}\log R_N}\ \hphi\left( \frac{\log p}{\log R_N}\right)
\nonumber\\ & = & r_\F \sum_{p=2}^{R_N^\sigma} \frac{\log p}{p\log
R_N}\ \hphi\left( \frac{\log p}{\log R_N}\right) +
O\left(\frac1{\afn^{\delta_1}}\sum_{p=2}^{R_N^\sigma} p^{\mu_1-\foh}
\right) \nonumber\\ & = & r_\F \phi(0) + O\left(\frac{1}{\log
R_N}\right) +
O\left(\frac{R_N^{(\mu_1+\foh)\sigma}}{\afn^{\delta_1}}\right), \eea
where the main term in the last line is an immediate consequence of
the Prime Number Theorem (see Theorem \ref{thmprimesums} for a
proof); as $\afn \ge R_N^{\delta_0}$, for $\sigma$ sufficient small
(in terms of $\mu_1, \delta_1$ and $\delta_0$), the last error term
is negligible.

We are left with the contribution from the squares of the primes
(the $\nu = 2$ terms). As $\sum_{f\in\fn} b_f(p^2) = c_\f\afn +
O(\afn^{1-\delta_2}p^{\mu_2})$, for sufficiently small support, up
to a negligible term by Theorem \ref{thmprimesums} the resulting sum
over primes is $\frac{\phi(0)}2$. Thus the 1-level density satisfies
\be  D_{1,\mathcal{F}}(\phi) \ = \ \hphi(0) - c_\f \cdot
\foh \phi(0) + r_{\mathcal{F}}\phi(0), \ee which for small support
agrees with the $1$-level densities of \eqref{eq:1ldwithtestfn}
(trivially modified if there are $r_{\mathcal{F}}$ forced
eigenvalues at $1$).
\end{proof}

\begin{rek}[Support of the test functions]
  The allowable support of $\widehat\phi$ is determinable from the
  constants $\delta_i, \mu_i$. In general, the support will not be
  large enough to distinguish the three orthogonal densities, though
  it will suffice to distinguish unitary from symplectic from
  orthogonal.
\end{rek}

\begin{rek}[General $n$-level density]
It is natural to investigate the 2-level density to distinguish the
orthogonal groups. To do so requires two additional pieces of
information: (1) the distribution of signs of functional equations
in the family; (2) being able to average over the family
$b_f(p_1)b_f(p_2^2)$ and $b_f(p_1^2)b_f(p_2^2)$. The presence of
cross terms can seriously complicate matters, though fortunately in
all families considered to date these terms can be converted to
products of averages of single terms. For cuspidal newforms this
follows from the Petersson formula; for one-parameter families of
elliptic curves, if $p_1, \dots, p_k$ are distinct primes and
$r_1,\dots,r_k$ are integers, simple counting (see \cite{Mil2})
shows that
\begin{eqnarray}\label{eqeqone}
\sum_{t \bmod p_1\cdots p_k} a_{t}^{r_1}(p_1) \cdots
a_{t}^{r_k}(p_k)\ =\ \prod_{i=1}^k \sum_{t \bmod p_i}
a_t^{r_i}(p_i),
\end{eqnarray} reducing the analysis to single product terms. In
general (though see Remark \ref{rek:varanalcond} below), if we know the
distribution of signs we can determine the $2$-level densities, as
all that matters is knowing the fraction of $L(s,f)$ with even or
odd sign; however, the story is markedly different for third- and
higher-level densities. There we need to know significantly more; we
need to know which of the $L(s,f)$ have odd functional equation, and
we need to execute sums over \emph{just those} $L(s,f)$. Thus, for
families of elliptic curves the third- and higher-level densities are
beyond current techniques (except for constant-sign families); see
\cite{Mil1} for more details.
\end{rek}

\begin{rek}[Variation in the analytic conductors]\label{rek:varanalcond}
  We concentrate on the $1$-level density in this paper. There are two
  ways to normalize the zeros of an $L$-function in a family
  $\mathcal{F}_N$: we can use $\log Q_f$ or $\log R_N$. For the
  $1$-level density, it does not matter which one is used in the
  normalization; however, for the higher $n$-level densities we
  need to use the explicit formula multiple times. While the
  normalization $\log R_N$ greatly simplifies the $1$-level
  computations (because all test functions are scaled equally), in the
  general case we are forced to evaluate sums of $\log Q_f$ against
  the coefficients of the $L$-functions.  With additional work, these
  sums can often be handled (see \cite{Mil2} for the case of
  one-parameter families of elliptic curves).
\end{rek}


\bigskip
\section{Proof of the Main Results}\label{sec:proofmainresults}

\subsection{Preliminaries}\ \\

Assuming $\f$ and $\g$ are $\ntg$ families, by Theorem
\ref{thm:1ntg} we can determine their 1-level densities, and
associate a classical compact group to the family (uniquely in the
case of unitary and symplectic symmetry; for orthogonal symmetries, for small support the $1$-level density cannot distinguish $\soe$ from $\so$ from $\soo$). Assuming a few additional conditions, we can determine the symmetry group of the Rankin-Selberg convolution of the two families $\f$ and $\g$.

The following lemma is key to the proofs of our main results.

\begin{lem}\label{lemftimesg} Assume $\pi_1,\pi_2$ are automorphic
  cuspidal representations of ${\rm GL}_n(\bfA_\bfQ)$, ${\rm
    GL}_m(\bfA_\bfQ)$, respectively, and further that the functorial
  lift $\pi_1\times\pi_2$ to $\GL_{mn}(\bfA_\bfQ)$ exists.  Assume
  that, for some prime $p$, both $\pi_{1,p}$ and $\pi_{2,p}$ are
  unramified, and their corresponding Satake parameters are
  $\{\alpha_{\pi_1}(i)\}_{1\leq i\leq n}$ and
  $\{\alpha_{\pi_2}(j)\}_{1\leq j\leq m}$.  Then
  \begin{equation}
    b_{\pi_1 \times \pi_2}(p^\nu)\ =\ b_{\pi_1}(p^\nu) \cdot b_{\pi_2}(p^\nu).
  \end{equation}
\end{lem}
\begin{proof} Let $b_{\pi_1}(p^\nu)$ and $b_{\pi_2}(p^\nu)$ be as in
  Definition~\ref{defntg}.
  By the local Langlands correspondence, the Satake parameters for
  $\pi_{1,p}\times\pi_{2,p}$ are
\begin{equation}
  \{\alpha_{\pi_1\times \pi_2}(k)\}_{k=1}^{nm}\ =\ \{\alpha_{\pi_1}(i)
  \cdot\alpha_{\pi_2}(j)\}_{{1 \le i \le n \atop 1 \le j \le m}},
\end{equation}
which gives
\begin{align}
 b_{\pi_1 \times\pi_2}(p^\nu) &  =  \sum_{k=1}^{nm}
\ga_{\pi_1\times\pi_2}(k)^\nu\nonumber\\ & =  \sum_{i=1}^n
\sum_{j=1}^m \ga_{\pi_1}(i)^\nu \cdot \ga_{\pi_2}(j)^\nu \nonumber\\ &  =
 \sum_{i=1}^n \ga_{\pi_1}(i)^\nu \cdot \sum_{j=1}^m \ga_{\pi_2}(j)^\nu
\nonumber\\ &  =  b_{\pi_1}(p^\nu) \cdot b_{\pi_2}(p^\nu).
\end{align}
\end{proof}
\begin{rek}
  In Lemma~\ref{lemftimesg} above the assumption that
  $\pi_1\times\pi_2$ is automorphic is not really necessary in order
  to define the Satake parameters at~$p$ nor the coefficients
  $b_{\pi_1\times \pi_2}(p^\nu)$ of a hypothetical $\pi_1\times\pi_2$.
  Recall that our viewpoint is that the $L$-function
  $L(s,\pi_1\times\pi_2)$ should be automorphic and primitive, which
  is the case if $\pi_1,\pi_2$ are cuspidal, \emph{except} when
  $\pi_1\simeq\widetilde\pi_2$ (assuming unitary and of trivial
  central character).
\end{rek}

If $\pi$ is an automorphic cuspidal representation of ${\rm GL}_n$,
then $L(s,\pi)$ is primitive if and only if $\pi$ is primitive. For
two unitary cuspidal automorphic representations $\pi$ and $\pi'$ of
${\rm GL}_n$ with trivial central character, their convolution
$L(s,\pi \times \pi')$ is primitive if and only if $\pi'$ is
\emph{not} the contragredient of $\pi$. This is equivalent to the lift
$\pi\times \pi'$ being a cuspidal representation of ${\rm GL}_{n^2}$.

We now prove our main results.

\bigskip
\subsection{Convolving two families}\ \\

\begin{proof}[Proof of Theorem \ref{thm:ntgRSconv}] We first prove the theorem under the assumption that all convolved $L$-functions are primitive, and then handle the general case.

As we are assuming the convolved family is $\ntg$, in the new family $\fn \times \gm$  the conductors are
essentially constant, say $\log Q_{f\times g} = \log R_{N,M} +
o(\log R_{N,M})$. By the multiplicativity assumptions, it will be relatively
easy to evaluate \bea & & D_{1,\fn\times\gm}(\phi)\nonumber\\ &= &
\hphi(0)-2 \cdot \frac1{\afn\cdot\agm} \sum_{f\times
g\in\fn\times\gm} \sum_{\nu=1}^\infty \sum_{p=2}^{R_{N,M}^\sigma}
\frac{b_{f \times g}(p^\nu) \log p}{p^{\nu/2}\log
R_{N,M}}\ \hphi\left( \nu\frac{\log p}{\log R_{N,M}}\right) + o(1). \nonumber\\
\eea Some care is required for the $\nu \ge 3$ terms. We need to
show that these give a negligible contribution. If condition (4.iii)
holds for either of the $\ntg$ families (namely, that if we sum over
the family, we have a power savings in the family cardinality), then
this follows immediately. If not, we need $\delta_3(\f), r_3(\f)$
(and similarly for the family $\g$) to be such that summing the $\nu
\ge 3$ terms is negligible. This is always the case if we assume the
Ramanujan conjecture, condition (4.i).

We must determine the contributions from the $\nu = 1, 2$ terms. As
\be b_{f \times g}(p^\nu) \ = \ b_f(p^\nu) \cdot b_g(p^\nu) \ee
(Lemma \ref{lemftimesg}), we can execute the summations over $f \in
\fn$ and $g\in\gm$. The main term from $\nu = 2$ is
\be\label{eq:maintermtwistingtwontgoodfamilies} -2 \cdot
\frac{1}{\afn \cdot \agm} \sum_{p=2}^{R_{N,M}^\sigma} \frac{c_\f\afn
\cdot c_\g\agm \log p}{p\log R_{N,M}}\ \hphi\left( 2\frac{\log
p}{\log R_{N,M}}\right). \ee If the zeros were normalized by $\log
Q_{f\times g}$ instead of $\log R_{N,M}$, it would not be as easy to
compute the contributions because the Schwartz functions would be
evaluated at points depending on $Q_{f\times g}$. This is the main
reason we choose to normalize all zeros in a family by the same
quantity.

By the Prime Number Theorem (Theorem \ref{thmprimesums}), the main term of the sum in
\eqref{eq:maintermtwistingtwontgoodfamilies} equals \be - \frac{c_\f
\cdot c_\g}2 \ \phi(0), \ee and the error term is negligible. There are three other terms which
contribute in the $\nu=2$ case: $$\afn^{1-\delta_2(\f)} \cdot\agm
\cdot p^{\mu_2(\f)}, \ \  \afn\cdot \agm^{1-\delta_2(\g)}\cdot
p^{\mu_2(\g)}, \ \ \afn^{1-\delta_2(\f)}\cdot
\agm^{1-\delta_2(\g)}\cdot p^{\mu_2(\f)+\mu_2(\g)}.$$ As we divide
by $\afn \cdot \agm$, each of the three terms leads to a negligible
contribution for test functions with suitably small support.

We are left with handling the $\nu=1$ terms. If $r_\F$ or $r_\G = 0$
then we immediately see this term does not contribute for suitably
small support. For notational convenience we assume
\eqref{eq:onelevelsumbfp} and not \eqref{eq:onelevelsumbfpb} holds, as
the argument in each case is similar. We have \bea \sum_{f\times g \in
  \fn\times\gm} b_{f\times g}(p) & \ = \ & \left[ \sum_{f\in\fn}
  b_f(p)\right]
\cdot \left[  \sum_{g\in\gm} b_g(p)\right] \nonumber\\
& = & \frac{r_\F\cdot r_\G}{p}\afn\cdot\agm +
O\left(\afn^{1-\delta_{1,\F}}\cdot\agm^{1-\delta_{1,\G}}
  p^{\mu_{1,\F}+\mu_{1,\G}} \right)\nonumber\\ & & \ +\
O\left(\afn^{1-\delta_{1,\F}}\cdot\agm\cdot p^{\mu_{1,\F}-\foh}
\right) \nonumber\\ & & \ +\
O\left(\afn\cdot\agm^{1-\delta_{1,\G}}\cdot p^{\mu_{1,\G}-\foh}
\right).  \eea Summing over $p$, for test functions with small support
the three error terms do not contribute. The main term leads to \be -2
\sum_p \frac{r_\F \cdot r_\G}{p} \frac{\log p}{\sqrt{p}\log R_{N,M}}\
\hphi\left(\frac{\log p}{\log R_{N,M}}\right). \ee If $\F\times\G$
were to have rank, this sum would have to contribute. Comparing to
equation~\eqref{eq:onelevelsumbfpb}, the difference is that in the sum
above we have $\frac{r_\F\cdot r_\G}{p}$ instead of something like
$\frac{r}{\sqrt{p}}$ times $\frac{\log p}{\sqrt{p}\log R}\
\hphi\left(\frac{\log p}{\log R}\right)$. The presence of $p$ rather
than $\sqrt{p}$ in the denominator means this sum is of size $(\log
R_{N,M})^{-1}$ rather than of size $1$. This leads to a lower order
correction term to the $1$-level density of size $\frac{r_\F\cdot
  r_\G}{\log R_{N,M}}$.

We now remove the assumption that all the convolutions are primitive; i.e., we now allow a
  contragredient of an $f\in\f$ to be in $\g$. This can only happen if $m=n$. All we require is some control on the number
  such pairs $(f,\tilde f)$ and their contribution. As the families are $\ntg$, the multiplicities are bounded: $\mu_f \le \mu_\f$ and $\mu_g \le \mu_\g$. Thus the number of such pairs is trivially bounded by $\min(\mu_\g |\f_N|, \mu_\f |\g_M|) = O(\min(|\f_N|,|\g_M|))$.

  For any such $(f,\tilde f)$ the
  convolution $L(s,f \times \tilde f)$ is not primitive and has a simple pole at $s=1$, contributing two additional terms, $\phi\left(\pm \frac{\log R}{4\pi}i\right)$, to the explicit formula (see, for example, \cite{RS}). If $\supp(\hphi) \subset (-\sigma, \sigma)$, then \bea \phi(t+iy) & \ = \ & \int_{-\infty}^\infty \hphi(\xi)
e^{2\pi i (t+iy)\xi} d\xi \ \ll \ e^{2\pi |y| \sigma}. \eea Thus \be \phi\left(\pm \frac{\log R}{4\pi}i\right) \ \ll \ R^{\sigma/2}; \ee as we divide by $|\f_N| \cdot |\g_M|$ and there are only $O(\min(|\f_N|,|\g_M|))$ such pairs, for $\sigma$ sufficiently small these two terms have a negligible contribution.

  We now show the contribution to the prime sums from these pairs is also negligible if $\sigma$ is sufficiently small. Any improvement of the exponent $1/2$ in the Jacquet-Shalika \cite{JS} bound for the Satake parameters suffices; we use the Rudnick-Sarnak\footnote{The Rudnick-Sarnak bound is stated only for Satake parameters of cuspidal representations, and it trivially extends to isobaric sums of cuspidal representations, hence to arbitrary (not necessarily cuspidal) automorphic representations. } \cite{RS}   bound: if $\pi$ is an automorphic representation of ${\rm GL}_r(\A_\Q)$, then $|\alpha_{\pi,j}(p)| \le p^{\foh - \frac1{r^2+1}}$. Thus, each pair contributes at most \be \sum_\nu \sum_{p \le R^{\sigma/\nu}} \frac{p^{\frac{\nu}2-\frac{\nu}{n^4+1}}}{p^{\nu/2}} \ \ll_n \ R^\sigma. \ee
   to the prime sums in the explicit formula. As there are only $O(\min(|\f_N|,|\g_M|))$ pairs, for $\sigma$ sufficiently small these lead to negligible contributions upon dividing by the family's cardinality, $|\f_N| \cdot |\g_M|$.\end{proof}

 \begin{rek} The universality in Theorem \ref{thm:ntgRSconv} can be surprising
   at first. In determining the underlying classical compact group of the
   convolution of two families, all that matters are the distribution of signs
   of functional equations, the rank of the family, and the family averages of
   the $b_f(p)$'s and $b_g(p)$'s (i.e., the family averages of the second
   moments of the Satake parameters at each unramified prime).  Upon convolving
   two such nice families, the main term is independent of the family ranks;
   however, there \emph{is} a lower order correction term which can often be
   isolated and which \emph{does} depend on the ranks. Unfortunately the bounds
   for the errors from the $\nu \ge 3$ terms, even assuming Ramanujan, will be
   of the same size, as could the other error bounds from the $\nu=1$ and
   $\nu=2$ terms. Conjecturally, however, it is reasonable to expect there to be
   cancelation in these errors upon summing over the families, and hence that
   there could be corrections to the $1$-level density. For other examples of
   lower order corrections, see \cite{FI,Mil3,Mil5,Mil6,Mil7,St,Yo1}. \end{rek}

\begin{exa}\label{exa:holocusp}
  We give an interesting example of Theorem \ref{thm:ntgRSconv}.
  Consider families $\F_i$ of weight-$k_i$ holomorphic cuspidal
  newforms of prime level~$N$ ($k_i$ fixed, $N \to \infty$); perhaps
  we might want to take the sub-families of even or odd sign. These
  families $\F_i$ are $\ntg$, as is $\F_1\times\F_2$ when $k_1 \neq
  k_2$. By \cite{ILS} each $\F_i$ has orthogonal symmetry. As these
  are $\glt$ holomorphic cuspidal newforms, we know the Ramanujan
  conjectures and thus condition (4.i) holds. From the Petersson
  formula, \eqref{eq:onelevelsumbfp} holds with $r_{\F_i} = 0$. Thus
  these families have orthogonal symmetries and hence their symmetry
  constants are $c_{\F_i} = -1$. Therefore $c_{\F_1\times\F_2} =
  c_{\f_1} \cdot c_{\f_2} = 1$, implying that $\F_1\times\F_2$ has
  symplectic symmetry. In particular, all elements should have even
  sign (which we do get from Rankin-Selberg). Note this is a
  $\mbox{GL}_4$ family of $L$-functions. Is there a larger natural
  $\mbox{GL}_4$ family containing it (analogous to the quadratic
  Dirichlet characters sitting inside all Dirichlet characters)?
  Further, if $k_1,k_2$ and $k_3$ are distinct then
  $c_{\F_1\times\F_2\times\F_3} = c_{\f_1} \cdot c_{\f_2} \cdot
  c_{\f_3} = -1$, implying $\F_1\times\F_2\times\F_3$ has orthogonal
  symmetry.\footnote{It is known that $f_1 \times f_2$ is automorphic.
    While it is not known that $f_1 \times f_2 \times f_3$ is
    automorphic, we do know that $L(s,f_1 \times f_2 \times f_3)$ is
    entire; the automorphicity follows from standard functoriality
    conjectures, and would imply the $L$-function is primitive. See
    \cite{Bu,Ga,Ram} for details.}
\end{exa}

\bigskip
\subsection{Convolving by a Fixed Form}\label{sec:twistfixed}\ \\

\begin{proof}[Proof of Theorem \ref{thm:twistfixedform}] We may assume all the convolutions are primitive. As $\g$ is $\ntg$, from our cardinality assumption there are at most $\mu_G = O(1)$ imprimitive convolutions. Arguing as in the proof of Theorem \ref{thm:ntgRSconv} we can show they have a negligible contribution for sufficiently small support.

Using the explicit formula to calculate the 1-level density, we
have the expansion \bea D_{1,f\times\gm}(\phi) &\ = \ & \hphi(0)-
\frac2{\agm} \sum_{g\in\gm} \sum_{\nu=1}^\infty
\sum_{p=2}^{R_M^\sigma} \frac{b_f(p^\nu)b_g(p^\nu) \log
p}{p^{\nu/2}\log
R_M}\ \hphi\left( \nu\frac{\log p}{\log R_M}\right) + o(1). \nonumber\\
\eea There will be no contribution from the $\nu \ge 3$ terms if we
have sufficiently good bounds for $\frac{b_f(p^\nu)
b_g(p^\nu)}{p^{\nu/2}}$, or if we have some power savings (relative
to $\agm$) in $\sum_{g\in\gm} b_g(p^\nu)$. For example, if we take
$f$ to be any nice $L$-function on $\mbox{GL}_2$ (say a holomorphic
cuspidal newform of weight $k$ and level $N$ or an even Maass form),
then we have good bounds on $b_f(p^\nu)$. In the holomorphic case,
we know Ramanujan and $b_f(p^\nu) \ll 1$; in the Maass case we have
$b_f(p^\nu) \ll p^{7/64}$ (see the appendix by Kim and Sarnak in \cite{K}). We can quantify exactly
what bounds we need on $b_g(p^\nu)$ for each $g\in\gm$ ($\nu \ge
3$), and these bounds are available in many cases of interest. We
then execute the summation over $p$, which gives $(\log R_M)^{-1}$,
and then we trivially handle the sum over $\gm$.

As $\g$ satisfies \eqref{eq:onelevelsumbfp}, for $\nu = 1$ a simple
calculation shows there is no contribution in the limit as $M \to
\infty$ for sufficiently small support. We are left with the crucial
case of $\nu = 2$; note that this is the term that determines the
symmetry type: If it is $0$ (resp., $1$ or $-1$), we have unitary
(resp., symplectic or orthogonal). From our assumption that $\g$
satisfies \eqref{eq:3iione}, we can execute the summation
$\sum_{g\in\gm} b_g(p^2)$, and we find that the main contribution from
$\nu = 2$ is just \be -\frac2{\agm} \sum_{p=2}^{R_M^\sigma}
\frac{\left[c_\g \cdot \agm + \agm^{1-\delta_2}
    p^{\mu_2}\right]b_f(p^2) \log p}{p\log R_M}\ \hphi\left(
  2\frac{\log p}{\log R_M}\right). \ee For sufficiently small support,
as $b_f(p^2)$ is bounded by some power of $p$, the second term doesn't
contribute. We are left with \be -2 c_\g \sum_{p=2}^{R_M^\sigma}
\frac{b_f(p^2) \log p}{p\log R_M}\ \hphi\left( 2\frac{\log p}{\log
    R_M}\right). \ee

Thus the symmetry will be the product of $c_\g$ and the above sum.
If $f$ is a Dirichlet character $\chi$, then $b_f(p^2) = \chi(p)^2$.
If $\chi$ is quadratic than $\chi(p)^2 = 1$ and the symmetry
constant will be $c_\g$ again; if $\chi$ is not quadratic than the
sum of $\chi(p)^2$ (times the other factors) over the primes is
$o(1)$, yielding unitary symmetry. If $f$ is a nice $\mbox{GL}_2$
$L$-function (say holomorphic cuspidal Hecke newform or
Hecke-Maass), then the prime sum is $-\fof\phi(0)$ because
\bea\label{eq:goodeqbfp2} b_f(p^2) & \ = \ & \ga_{f,1}^2(p) +
\ga_{f,2}^2(p) \nonumber\\ & = & (\ga_{f,1}(p) + \ga_{f,2}(p))^2 - 2
\nonumber\\ & = & a_f(p)^2 - 2 \nonumber\\ & = & \left[a_f(p^2) +1
\right] - 2 \nonumber\\ & = & a_f(p^2) - 1, \eea where we have used
the fact that $f$ is a Hecke eigenform to say $a_f(p)a_f(p) =
a_f(p^2) + 1$ (at least for $p$ relatively prime to the conductor).
The $a_f(p^2)$ will be related to the symmetric square $L$-function
associated to $f$, and by GRH for that $L$-function, its sum over
primes is negligible (see \cite{ILS} for details). Thus the $\nu =
2$ terms contribute $(-2c_\g) \cdot (-\fof \phi(0)) = c_\g \cdot
\foh \phi(0)$.

Setting $c_f = 1$ if $f$ is a quadratic Dirichlet character, $0$ if
$f$ is a non-quadratic Dirichlet character, and $c_f = -1$ if $f$ is
a Hecke holomorphic or Maass form, we find that the $1$-level
density of $f\times\g$ is \be \hphi(0) - c_f \cdot c_\g \cdot \foh
\phi(0). \ee  \end{proof}


\begin{rek} These results are similar to those obtained by Rubinstein
  in his thesis, where he considered the convolution of the family of
  quadratic Dirichlet $L$-functions with a fixed ${\rm GL_n}$ form;
  see \cite{Rub}. In our notation, if $f$ is self-dual, then $c_f=+1$
  (resp., $c_f=-1$) if $L(s,\sym^2f)$ (resp., $L(s,\wedge^2f)$) has a
  pole at $s=1$. If $f$ is not self-dual then $c_f=0$.
\end{rek}


\bigskip
\section{Convolving Families Of Symmetric Powers Of Modular Forms}\label{sec:symmpowers}

Families of $L$-functions attached to holomorphic modular forms and their
functorial liftings are often $\ntg$, at least under the assumption of standard
conjectures.  The main purpose of this section is to provide further examples
illustrating Theorem~\ref{thm:ntgRSconv}. Additional examples (independent of
this section) involving elliptic curves are given in \S\ref{sec:ellcurves}.

Let $\cH_k$ be a Hecke eigenbasis of the space of modular cusp forms
of weight~$k$ for the full modular group $\text{SL}_2(\bfZ)$.  Then
$|\cH_k|=\frac k{12}+O(1)$. We denote the average over $\cH_k$ by
\be \langle A_f \rangle_{\cH_k} \ := \ \frac1{|\cH_k|} \sum_{f\in
\cH_k} A_f. \ee

We normalize $f\in\cH_k$ so its leading Fourier coefficient is one,
viz.,
\begin{align}
  \label{eq:Fourier_L_modular}
  f(z) & \ = \  \sum_{n=1}^\infty a_f(n)n^{\frac{k-1}2}\exp(2\pi inz) \\
  L(s,f) & \ = \  \sum_{n=1}^\infty a_f(n)n^{-s} \ = \  \prod_p(1-a_f(p)
  p^{-s}+p^{-2s})^{-1}\\
  & \ = \
  \prod_p(1-\alpha_f(p)p^{-s})^{-1}(1-\alpha_f(p)^{-1}p^{-s})^{-1},
  \qquad\Re s>1, \\
  \label{eq:FE-modular}
  \Lambda(s,f) & \ = \  2(2\pi)^{-\left(s+\frac{k-1}2\right)}
  \Gamma\left(s+\frac{k-1}2\right) L(s,f) \ = \  (-1)^{k/2}\Lambda(1-s,f).
\end{align}
Here $\alpha_f(p),\alpha_f(p)^{-1}$ are the Satake parameters
at~$p$. Since we never need to look at Satake parameters
simultaneously for two different primes, we usually omit $p$ and write
simply $\alpha_f$.  It is well known that $f$ uniquely determines an
automorphic cuspidal unitary \emph{self-dual} representation $\pi$ of
$\glt$ with trivial central character~\cite{Gel}.  Moreover $\pi_\infty$ is the
discrete series representation of weight~$k$.\footnote{Some authors
prefer to say that this $\pi_\infty$ has weight $k-1$.  We follow the
convention in~\cite{CM}.} In what follows we will implicitly use this
identification and rarely bother to talk about the representation
$\pi$ \emph{per se.}  From the completed $L$-function in
equation~\eqref{eq:FE-modular} (in particular from its gamma factor)
it follows that the analytic conductor of $\mathcal{F}_k$ is
$R_k\asymp k^2$.

Because $\cH_k$ consists of forms of full level, $\pi_p$ is
unramified for all~$p$.  The following orthogonality relations for
the Fourier coefficients $\{a_f(n)\}$ are crucial:
\begin{lem} We have
  \begin{equation}
  \label{eq:Petersson}
  \frac1{|\cH_k|+O(1)}\sum_{f\in\cH_k}\frac{\zeta(2)}{L(1,\sym^2f)}\ a_f(m)a_f(n)
  \ =\ \delta(m,n) + \mathcal E,
\end{equation}
where
\begin{equation}
  \label{eq:delta}
  \delta(m,n)\ =\
  \begin{cases}
    1& m=n\\
    0& m\neq n
  \end{cases}
\end{equation}
and
\begin{equation}
  \label{eq:Petersson_error}
  \mathcal{E}\ =\
  \begin{cases}
    O_\ell\left(\frac{(mn)^{1/4}\log mn}{k^{5/6}}\right)& \text{{\rm if $m,n$
    have no more than $\ell$ factors}}\\
  O\left(\frac{\sqrt{mn}}{2^k}\right)&\text{{\rm if} $12\pi\sqrt{mn}\leq k$.}
  \end{cases}
\end{equation}
\end{lem}
Formula \eqref{eq:Petersson} is a consequence of the Petersson
formula (see equation 2.12 of \cite{ILS}).  Note that the left-hand
side of \eqref{eq:Petersson} is just the average $\langle
a_f(m)a_f(n)\rangle_{\cH_k}$, except for the presence of the weights
$\zeta(2)/L(s,\sym^2f)$.  This is called the \emph{harmonic
averaging} of $a_f(m)a_f(n)$ and often makes the analysis more
tractable (see~\cite{DM,ILS,Mil7,Ro}). If we were interested in bounding
the order of vanishing at the central point in the family then the
harmonic weights would cause difficulty (see Remarks 2.11 and 6.1
in~\cite{HM}).

Following \cite{ILS}, by additional work we can remove the harmonic
weights in the $1$-level density. The cost is a slight worsening of the
constants $\delta_1,\delta_2,\delta_3$ in the definition of
$\ntg$. Alternatively, we can simply redefine the average $\langle
a_f(m)a_f(m)\rangle_{\cH_k}$ to be given by the left-hand side of
\eqref{eq:Petersson}.

Note that \be a_f(p^n)\ = \
\alpha_f^n+\alpha_f^{n-2}+\dots+\alpha_f^{-n+2}+\alpha_f^{-n}, \ee
so from Definition~\ref{def:b} it follows immediately that
\begin{align}
  \label{eq:b_modr_vs_a_modr}
  b_f(p) &\ =\ a_f(p)\nonumber\\
  b_f(p^2) &\ = \ a_f(p^2)-1.
\end{align}

These formulas, together with the orthogonality
relations~(\ref{eq:Petersson}), already suffice to prove conditions
(3.i) and (3.ii) of Definition~\ref{defntg} with
$\delta_1=\delta_2=1/6$, any $\mu_1>1/4$, $\mu_2>1/2$, rank zero and,
most importantly, with symmetry constant~$-1$ (note that
$a_f(p^n)=a_f(p^n)a_f(1)$ and $\delta(1,p^n)=0$ for $n=1,2$, whereas
$-1=-a_f(1)a_f(1)$).  Conditions (1) and~(2) are obvious, and the
Ramanujan conjecture (condition (4)) is known for these $f$ by
Deligne. We therefore recover the result from \cite{ILS,Ro} that the
family $\{\cH_k\}$ as $k\to\infty$ has orthogonal $1$-level density
(at least for small support).

For small support of test functions, one cannot in general pinpoint
the exact underlying symmetry in the orthogonal case. However, with
the help of the root number (sign of the functional equation), the
symmetry should be SO(even) if all the functional equations have
positive sign and SO(odd) if all have negative sign. Determining the
sign of the functional equation is most easily done through the
local Langlands correspondence. Since we will be building
automorphic representations starting from modular forms of full
level, all finite places ($p$ prime) contribute local root numbers
equal to~$+1$, and we only need the archimedean local
correspondence. Moreover, since the only archimedean place of $\bfQ$
is $\bfQ_\infty=\bfR$ we can simplify the notation a bit. The reader
who wants an authoritative survey of the archimedean Langlands
correspondence should read Knapp's article~\cite{Kn}.

The archimedean local correspondence for $\gln(\bfA_\bfQ)$ is a
bijection $\rho\leftrightarrow\pi_\infty$ between admissible
representations $\rho:\WR\to\gln(\bfC)$ and
\emph{irreducible} admissible representations $\pi_\infty$ of
$\gln(\bfR)$.  Here $\WR:=\bfC^\times\cup j\bfC^\times$
(disjoint union) is a multiplicative group with $j^2=-1$, and $j$ acts
on $\bfC^\times$ by $jzj=\bar z$.  $\WR$ is the \emph{Weil
  group} of $\bfR$; it can also be identified with an obvious
multiplicative subgroup of the quaternions.  We will not discuss the
meaning of admissibility here.

Irreducible admissible representations of $\WR$ are one or two
dimensional.  There are two families of inequivalent one-dimensional
representations, each parametrized by a complex number $t\in\bfC$.
They are denoted $\{[+,t]\}$ and $\{[-,t]\}$.  Additionally, there are
two-dimensional representations; they are parametrized by an integer
$k\geq2$ and a complex number $t\in\bfC$.  They are denoted $[k,t]$.
There are no irreducible admissible representations of dimension
greater than two, and any (finite-dimensional) admissible
representation of $\WR$ is fully reducible (decomposes as a direct sum
of irreducible ones).

The correspondence assigns $[+,0]$ to the trivial representation and
$[-,0]$ to the ``sign'' representation
$x\mapsto\sgn(x)=x\left|x\right|^{-1}$ of $GL(1,\bfR)$.  The
discrete-series representation of weight $k\geq2$ corresponds to
$[k,0]$.  The parameter $t\in\bfC$ parametrizes twists: either by the
character $\left|x\right|^t$ of $GL(1,\bfR)$ or by
$\left|\det(x)\right|^t$ of $GL(2,\bfR)$.

In order to characterize the archimedean components of functorial
liftings of automorphic representations, we need to understand the
effect of certain operations on representations of $\WR$.
\begin{lem} \label{lem:ops_weil_repns} Let $(-)^\kappa$ be `$+$' for
  $\kappa$ even, `$-$' for $\kappa$ odd.  Then for all $m\geq1$,
  $k>k'\geq2$, $t,t'\in\bfC$:
    \begin{align}
      \wedge^2[k,t] &\simeq [(-)^k,2t]\\
      \sym^m[+,t] &\simeq [+,mt]\\
      \sym^m[-,t] &\simeq [(-)^m,mt]\\
      \label{eq:sym-odd}
      \sym^{2m+1}[k,t] & \simeq \bigoplus_{\ell=0}^m
      [(2\ell+1)(k-1)+1,(2m+1)t]\\
      \label{eq:sym-even}
      \sym^{2m}[k,t] & \simeq
      [(-)^{m(k-1)},2mt]\oplus\bigoplus_{\ell=1}^m[2\ell(k-1)+1,2mt]\\
      [+,t]\otimes[+,t'] &\simeq [-,t]\otimes[-,t] \simeq [+,t+t']\\
      [+,t]\otimes[-,t'] &\simeq [-,t']\otimes[+,t] \simeq [-,t+t']\\
      [+,t]\otimes[k,t'] &\simeq [-,t]\otimes[k,t'] \simeq [k,t+t']\\
      \label{eq:convol-different}
      [k,t]\otimes[k',u] &\simeq [k',u]\otimes[k,t] \simeq
      [k+k'-1,t+t']\oplus[k-k'+1,t+t'] \\
      \label{eq:convol-equal}
      [k,t]\otimes[k,t'] &\simeq [2k-1,t+t']\oplus[+,t+t']\oplus[-,t+t']
    \end{align}
\end{lem}
The proof is easy and we omit it.  Cogdell and Michel
prove~(\ref{eq:sym-odd}) and~(\ref{eq:sym-even}) in~\cite{CM}.

The archimedean $\varepsilon$- and $L$- (gamma) factors are as
follows\footnote{With respect to standard Lebesgue measure on
$\bfR$ and provided the additive character used to define the
Fourier transform is $x\mapsto e^{2\pi ix}$.} (see \cite{Kn}):
\begin{align}
  \label{eq:Gamma-R}
  \GR(s) &:\ = \  \pi^{-s/2}\Gamma(s/2)\\
  \label{eq:Gamma-C}
  \GC(s) &:\ = \  \GR(s)\GR(s+1) \ = \  2(2\pi)^{-s}\Gamma(s) \\
\label{eq:L_epsilon_+}
  L(s,[+,t]) &\ = \  \Gamma_\bfR(s+t) & \varepsilon([+,t]) &\ = \  1 \\
\label{eq:L_epsilon_-}
  L(s,[-,t]) &\ = \  \Gamma_\bfR(s+t+1) & \varepsilon([-,t]) &\ = \  i \\
\label{eq:L_epsilon_k}
  L(s,[k,t]) &\ = \  \Gamma_\bfC\left(s+t+{\textstyle\frac{k-1}2}\right)
  & \varepsilon([k,t]) &\ = \  i^k.
\end{align}

Finally, $\varepsilon$- and $L$-factors are multiplicative with
respect to direct sums of representations of $\WR$ (which, via the
archimedean Langlands correspondence, are associated to isobaric
sums of irreducible admissible representations of
$\GL_{n_i}(\bfR)$), and if $\rho\leftrightarrow\pi_\infty$ then
$L(s,\pi_\infty)=L(s,\rho)$, and similarly for
$\varepsilon$-factors.

With these results in hand we can easily determine the underlying
symmetry type of various families obtained by functorial operations
starting from $\{\cH_k\}$.  However, we introduce one last bit of
notation: since, for $f\in\cH_k$, the automorphic representation
$\pi_f$ is unitary and has trivial central character, the parameter
$t\in\bfC$ is always zero in our applications and we will adopt the
following:\\

\textbf{Convention:} We write $[k]$ for $[k,0]$, $[+]$ for $[+,0]$,
and $[-]$ for $[-,0]$.  Also define $[1,t]:=[+,t]\oplus[-,t]$ and
$[1]:=[1,0]$.  Then equation~\eqref{eq:convol-equal} is the special
case $k=k'$ of~\eqref{eq:convol-different}.  Note that $[1,t]$ is a
\emph{reducible} two-dimensional representation of $\WR$.\\

\subsection{Families of Symmetric Powers}
\label{sec:fam-sym-powers}

We begin with the family $\mathcal{G}^{(M)}_k=\sym^M\cH_k$ for a
fixed $M\geq1$, and study the limit as $k\to\infty$.  It is
conjectured, and we assume this as a hypothesis, that every
$f\in\cH_k$ has a self-dual automorphic cuspidal functorial lift
$g=\sym^Mf\simeq\otimes'(\sym^Mf_v)$ whose local factors $\sym^Mf_v$
are defined through the local Langlands correspondence (by
composition with the $M$-th symmetric-power of the defining
representation of $GL_2(\bfC)$).  This is known for $M=1,2,4$ by
work of Hecke, Gelbart-Jaquet, Kim-Shahidi, and
Kim~\cite{GJ,KiSh1,K}.  Under this hypothesis, the family
$\mathcal{G}^{(M)}_k=\sym^M\cH_k$ consists of primitive
$L$-functions for $GL_{M+1}$.

\begin{thm}
  With the assumptions above, the family $\mathcal{G}^{(M)}$ is
  $\ntg$ with symmetry constant $c_{\mathcal{G}}=(-1)^M$.
\end{thm}
\begin{proof}
  Firstly, $|\mathcal{G}_k^{(M)}|=|\cH_k|=k+O(1)\to\infty$ as
  $k\to\infty$, so the cardinality condition holds.

  Now, $f_\infty\simeq[k]$ for all $f\in\cH_k$, so
  $g_\infty\simeq\sym^Mf_\infty$ are all \emph{isomorphic} admissible
  representations of $\GL_{M+1}(\bfR)$ as $f$ varies over $\cH_k$.  In
  addition, all the non-archimedean places are unramified, so the
  analytic conductors $Q_g$ are completely determined by $g_\infty$,
  and hence are \emph{constant} in $\mathcal{G}_k$.

  To compute the $\varepsilon$- and $\Gamma$-factor $L_\infty(s,\sym^M
  f)$ we use the Langlands correspondence,
  Lemma~\ref{lem:ops_weil_repns}, and equations
  \eqref{eq:L_epsilon_+}--\eqref{eq:L_epsilon_k}.  Recall that
  $f_\infty\simeq[k]$.

  We split into two cases: $M=2m$ and $M=2m+1$.
  \begin{itemize}
  \item \underline{$M=2m+1$}.
    \begin{align}
      \label{eq:gamma-sym-odd}
      \Gamma(s,\sym^{2m+1}f) &=
      \prod_{\ell=0}^m\GC\left(s+{\textstyle\frac{(2\ell+1)(k-1)}{2}}\right)\\
      \varepsilon(s,\sym^{2m+1}f) &=
      \begin{cases}
        i^k & m\equiv0\pmod4 \\
        -1 & m\equiv1\pmod4 \\
        -i^k & m\equiv2\pmod4 \\
        1 & m\equiv3\pmod 4.
      \end{cases}
    \end{align}
  \item \underline{$M=2m$}.
    \begin{align}
      \label{eq:gamma-sym-even}
      \Gamma(s,\sym^{2m}f) &=
      \GR\left(s+{\textstyle\frac{1-(-1)^{m(k-1)}}{2}}\right)
      \prod_{\ell=1}^m\GC\left(s+{\textstyle\frac{2\ell(k-1)}{2}}\right)\\
      \varepsilon(s,\sym^{2m}f) &= 1.
    \end{align}
  \end{itemize}
  As explained in~\cite{DM}, the contribution of the archimedean
  places to the analytic conductor $Q_g$ can be read off from the
  gamma factors (cf., equations~\eqref{eq:14} and~\eqref{eq:afqfrN} at
  the end of \S\ref{sec:central-1-2-level}): each factor $\GR(s+T)$
  contributes a factor of $T/2$ to the analytic conductor, and each
  factor $\GC(s+T)$ contributes $T(T+1)/4\asymp T^2/4$.
  Equations~\eqref{eq:gamma-sym-odd} and~\eqref{eq:gamma-sym-even}
  reveal that the analytic conductor is $Q_g\asymp(k/2)^{M+1}$ if $M$
  is odd, $Q_g\asymp(k/2)^{M}$ if $M$ is even.  This verifies the
  conductors condition.  The error terms are handled using the Ramanujan
  bounds of Deligne.

  To analyze the crucial conditions on prime sums, let us write
  $\alpha_p,\alpha_p^{-1}$ for the Satake parameters of $f_p$.  Then those
  of $(\sym^M f)_p$ are
  \begin{equation}
    \label{eq:Satake-sym-M}
    \alpha_p^M,\ \alpha_p^{M-2},\ \dots,\ \alpha_p^{-M+2},\ \alpha_p^{-M}.
  \end{equation}
  Writing
  $a(p^n),b(p^n)$ for $a_f(p^n),b_f(p^n)$ and $B(p^n)$ for
  $b_{\sym^Mf}(p^n)$, we have:
  \begin{align}
    B(p) &\ =\ \alpha_p^M + \alpha_p^{M-2}+\dots+\alpha_p^{-M+2} + \alpha_p^{-M}\notag\\
    \label{eq:B-sym-m-p}
    &\ =\ a(p^M) \\
    B(p^2) &\ =\ \alpha_p^{2M} + \alpha_p^{2M-4} + \dots + \alpha_p^{-2M+4} +
    \alpha_p^{-2M}\notag\\
    \label{eq:B-sym-m-p2}
    &\ =\ a(p^{2M})-a(p^{2M-2})+a(p^{2M-4})-a(p^{2M-6})+\dots+(-1)^M a(1).
  \end{align}
  Once again the orthogonality relations of~(\ref{eq:Petersson}) prove
  the condition on the prime sum (with rank zero) and the prime square sum,
  with symmetry constant $c_{\sym^M\cH_k}=(-1)^M$.  This reveals
  underlying symplectic symmetry when $M=2m$ is even and orthogonal
  when $M=2m+1$ is odd.  In the latter case, by looking at the
  $\varepsilon$-factor (root number), we expect that the symmetry is
  SO(even) for $m\equiv3\bmod4$, SO(odd) for $m\equiv1\bmod4$, and
  full orthogonal when $m$ is even. Furthermore, in this last case the
  symmetry is SO(even) (resp., SO(odd)) when $k/2$ is even (resp.,
  $k/2$ is odd).
\end{proof}
\begin{rek}
  G\"ulo\u{g}lu has obtained results for larger support for
  symmetric-power families~\cite{Gu}.
\end{rek}

\subsection{Convolutions of Symmetric Powers}
\label{sec:conv-sym-powers}

\begin{thm}\label{thm:sym_times_sym}
  Fix $M,N\geq1$ and consider the families
  $\mathcal{F}_k^{(M)}=\sym^M\cH_k$ and
  $\mathcal{G}^{(N)}_{k'}=\sym^N\cH_{k'}$.  Assume that the
  convolutions $f\times g$, $f\in\mathcal{F}_k^{(M)}$,
  $g\in\mathcal{G}_{k'}^{(N)}$ are automorphic.  Let
  $\mathcal{H}_{k,k'}^{(M,N)}=\mathcal{F}^{(M)}_k\times\mathcal{G}^{(N)}_{k'}$
  (where, as usual, we discard the non-cuspidal $f\times f$ when $M=N$
  and $k=k'$).  Then, as $k,k'\to\infty$ in such a way that $\log
  k'/\log k\to 1$, the family $\mathcal{H}^{(M,N)}$ $=$ $\mathcal{F}^{(M)}\times\mathcal{G}^{(N)}$
  is $\ntg$ with
  symmetry constant $c_{\mathcal{H}^{(M,N)}}=(-1)^{M+N} = c_{\mathcal{F}^{(M)}}
  \cdot c_{\mathcal{G}^{(N)}}$.
\end{thm}

\begin{rek}
  The automorphicity of the convolutions $f\times g$ is known when
  $M+N\leq 3$~\cite{Ram,KiSh1}.
\end{rek}

\begin{proof}
  For simplicity we will only consider the case when
  $k=k'\to\infty$; the proof of the general case differs from this case only in
  trivial details.

  As in the previous section, all non-archimedean places of
  $f\in\mathcal{F}^{(M)}_k$, $g\in\mathcal{G}^{(N)}_k$ are unramified.  We will
  once more split into cases when $M,N$ are even or odd.

  Using Lemma~\ref{lem:ops_weil_repns} we obtain%
  \footnote{Recall $[1]:=[+]\oplus[-]$, and observe that
    $[+]\otimes[1]\simeq[-]\otimes[1]\simeq[1]$.}:
  \begin{align}
    \label{eq:conv-symm-odd}
    \sym^{2m+1}[k]\otimes\sym^{2n+1}[k]&\simeq \bigoplus_{{0\leq\ell\leq
        m\atop 0\leq\lambda\leq n}} \bigl([2(\ell+\lambda+1)(k-1)+1]
    \oplus[2\left|\ell-\lambda\right|(k-1)+1]\bigr) \\
    \sym^{2m}[k]\otimes\sym^{2n}[k]&\simeq
    [(-)^{m+n}]\oplus\bigoplus_{1\leq\ell\leq m}[2\ell(k-1)+1]
    \oplus\bigoplus_{1\leq\lambda\leq n}[2\lambda(k-1)+1] \notag \\
    \label{eq:conv-symm-even}
    & \quad\oplus \bigoplus_{{1\leq\ell\leq m\atop 1\leq\lambda\leq
        n}}
    \bigl([2(\ell+\lambda)(k-1)+1]
    \oplus[2\left|\ell-\lambda\right|(k-1)+1]\bigr)
  \end{align}
  \begin{multline}
    \label{eq:conv-symm-even-odd}
    \sym^{2m+1}[k]\otimes\sym^{2n}[k]\simeq
    \bigoplus_{0\leq\ell\leq m}[(2\ell+1)(k-1)+1] \\
    \quad\oplus \bigoplus_{{0\leq\ell\leq m\atop 1\leq\lambda\leq
        n}}
    \bigl([(2\ell+2\lambda+1)(k-1)+1]
    \oplus[\left|2\ell-2\lambda+1\right|(k-1)+1]\bigr).
  \end{multline}
  The $\varepsilon$-factors are as follows:
  \begin{align}
    \label{eq:epsilon-RS-equal}
    \varepsilon(\sym^{2m+1}[k]\otimes\sym^{2n+1}[k]) &=
    \varepsilon(\sym^{2m}[k]\otimes\sym^{2n}[k])\ =\ +1 \\
    \label{eq:epsilon-RS-different}
    \varepsilon(\sym^{2m+1}[k]\otimes\sym^{2n}[k]) &\ =\
    \begin{cases}
      (-1)^{(m+1)(n-m)+(m+1)^2\frac k2} & m<n \\
      (-1)^{\frac{(m-n)(m+n+1)}2 + (m+1)^2\frac k2} & m\geq n.
    \end{cases}
  \end{align}


  We omit explicitly writing down the $\Gamma$-factors, but observe
  that every term $[a(k-1)+1]$ with $a>0$ contributes a factor $\asymp
  \frac{a^2}{4}\cdot k^2$ to the analytic conductor $Q_{f\times g}$.
  Hence, up to an additive constant, the analytic log-conductors are
  $\log Q_{f\times g}\sim 2(2m+1)(n+1)\log k$, resp.~$2m(2n+1)\log k$,
  resp.~$2(m+1)(2n+1)\log k$ corresponding to the
  cases~\eqref{eq:conv-symm-odd}, resp.~\eqref{eq:conv-symm-even},
  resp.~\eqref{eq:conv-symm-even-odd} above (we assumed $m\geq n$ in
  the first two cases).

  By Theorem~\ref{thm:ntgRSconv}, it only remains to show that
  $\mathcal{F}^{(M)}\times\mathcal{G}^{(N)}$ is $\ntg$.  The argument above
  shows that the conductor condition is satisfied.  When $M=N$ and
  $k=k'$ the representations $f\times f$ are not cuspidal; hence we
  must discard $O(k)$ of them.  Note that
  $|\mathcal{F}^{(M)}_k\times\mathcal{G}^{(N)}_k|=k^2+O(k)$ (so the cardinality
  condition holds) and that possibly shrinking the family introduces
  error terms of size $O(1/k)$, which are quite admissible.
  Properties of cardinality and the handling of error terms (by
  Ramanujan) are thus valid.  We need not verify the conditions on
  prime sums explicitly: the reason is that there are no ramified
  primes and conductors are essentially constant.  Hence
  Lemma~\ref{lemftimesg} and the argument in the proof of
  Theorem~\ref{thm:ntgRSconv} suffice to prove that
  $\mathcal{F}^{(M)}\times\mathcal{G}^{(N)}$ is $\ntg$ with symmetry
  constant $c_{\mathcal{F}^{(M)}\times\mathcal{G}^{(N)}}=c_{\mathcal{F}^{(M)}}\cdot
  c_{\mathcal{G}^{(N)}}=(-1)^{M+N}$.

  Therefore, for small support, the $1$-level density of the family
  $\mathcal{F}^{(M)}\times\mathcal{G}^{(N)}$ agrees with symplectic for $M+N$
  even, whereas for $M=2m+1$, $N=2n$ the symmetry is orthogonal and
  the root number, as read off from
  equations~\eqref{eq:epsilon-RS-equal}
  and~\eqref{eq:epsilon-RS-different}, determines whether the
  underlying symmetry is SO(even) or SO(odd).
  Equation~(\ref{eq:epsilon-RS-equal}) holds even when $k\neq k'$, but
  the form of equation~(\ref{eq:epsilon-RS-different}) is specific to
  the case $k=k'$.
\end{proof}




\bigskip
\section{Convolving Families Of Elliptic
Curves}\label{sec:ellcurves}

We now consider the interesting case of convolving two families of
elliptic curves. Specifically, consider the one-parameter families
\bea\label{eq:familiesFnGmec} \mathcal{F}_N: &\ y^2  \ = \ x^3 +
A_1(T)x + B_1(T), \ \ \ T \in [N,2N-1] \nonumber\\ \mathcal{G}_M:&\
y^2 \ = \  x^3 + A_2(S)x + B_2(S), \ \ \ S \in [M,2M-1], \eea where
the polynomials $A_1(T)$ through $B_2(S)$ have integer coefficients.
If we specialize $T$ to $t$ we obtain an elliptic curve
$E_{\mathcal{F}}(t)$ with discriminant $\Delta_{\mathcal{F}}(t)$ and
conductor $C_{\mathcal{F}}(t)$; similarly if we specialize $S$ to
$s$ we obtain an elliptic curve $E_{\mathcal{G}}(s)$ with
discriminant $\Delta_{\mathcal{G}}(s)$ and conductor
$C_{\mathcal{G}}(s)$. The conductors are products of powers of
primes dividing the discriminants. It is known (see
\cite{BCDT,TW,Wi}) that the $L$-function of an elliptic curve of
conductor $C$ agrees with a weight-$2$ cuspidal newform of level~$C$.
Thus if $E_i$ are elliptic curves with conductors $C_i$ and
associated newforms $f_i$, by $L(s,E_1 \times E_2)$ we mean the
Rankin-Selberg convolution $L(s,f_1\times f_2)$, which is a ${\rm
GL}_4$ $L$-function. The arithmetic conductor $Q(f_1\times f_2)$ of
such $L(s,f_1\times f_2)$ is an integer satisfying
\be\label{eq:sizecondsc1c2} (C_1C_2)^2 / (C_1,C_2)^4 \ \le \
Q(f_1\times f_2) \ \le \ (C_1C_2)^2 / (C_1,C_2), \ee where
$(C_1,C_2)$ is the greatest common divisor of $C_1$ and $C_2$; see
for example \cite{HaMi}. We often write $Q(C_1,C_2)$ for $Q(f_1
\times f_2)$.

We are interested in the behavior of $\fn \times \gm$ as $N$ and $M$
tend to infinity. The gamma factors for these ${\rm GL}_4$
$L$-functions depend neither on the specific curve nor on the
family. As such, since we need only identify the analytic conductor
up to a constant, we may use the integer $Q(C_1,C_2)$ as the
analytic conductor.

We normalize the low lying zeros for the convolution $L$-function by
the average of the logarithms of the analytic conductors. Thus, we
set \be\label{eq:avelogcondtwistec} \log R_{N,M} \ := \ \frac{1}{NM}
\sum_{t=N}^{2N-1} \sum_{s=M}^{2M-1} \log Q(C_\f(t), C_\g(s)). \ee

We need $R_{N,M}$ to tend to infinity with $N$ and $M$. A weak
estimate on the size of $Q(C_\f(t),C_\g(s))$, namely that the
average log-conductor in \eqref{eq:avelogcondtwistec} tends to
infinity with $N$ and $M$, suffices for our purposes.

To show this requires a few basic facts about elliptic curves. An
elliptic curve $E: y^2 = x^3 + a_4 x + a_6$ has discriminant $\Delta
= -16(4 a_4^3 + 27 a_6^2)$ and $j$-invariant $j = 3 a_4^3 / (4 a_4^3
+ 27 a_6^2)$; it is also convenient to set $c_4 = -48 a_4$ and $c_6
= -864 a_6$. Let $R$ be the ring of integers for some local field
$K$; $K$ is a local field which is complete with respect to a
discrete valuation $v$. Let $\mathcal{M} = \{x \in K: v(x) > 0\}$ be
the maximal ideal of $R$, and let $k = R/\mathcal{M}$ be the residue
field. If $a_i \in R$ and $v(c_4) < 4$ or $v(c_6) < 6$, then the
equation for the elliptic curve is minimal with respect to the
valuation $v$.

\begin{thm}\label{thm:sizearithcondprod}
Notation as above, assume that there are non-constant monic integral
polynomials $f_1(x)$ and $g_1(x)$ such that $f_1(x)$ divides
$\Delta_\f(x)$ and $g_1(x)$ divides $\Delta_\g(x)$. To simplify the
analysis, assume $f_1(x)$ does not divide either $c_{\f,4}(x)$ or
$c_{\f,6}(x)$ (and similarly for $g_1(x)$). Define the average
log-conductor by \eqref{eq:avelogcondtwistec}. If $j_\f(T)$ and
$j_\g(S)$ are both non-constant, then for some $a > 0$
\be\label{eq:lowerboundsrnmcond} \frac{\log NM}{(\log \log
\min(N,M))^a} \ \ll_{\mathcal{F},\mathcal{G}} \ \log R_{N,M}  \
\ll_{\mathcal{F},\mathcal{G}} \ \log NM. \ee
\end{thm}


The proof follows from basic facts on solutions to Diophantine
equations and properties of elliptic curves, and is given in
Appendix \ref{sec:proofthmsizearithcondprod}.

The following observation ensures that, except for a negligible
fraction of the time, the $L$-functions in the convolved family are
good (i.e., primitive).

\begin{lem}\label{lem:nozetamosttime}
Assume $j_\f(T)$ and $j_\g(S)$ are non-constant. The Rankin-Selberg
convolution of $E_\f(t)$ and $E_\g(s)$ is imprimitive for at most $O(\min(N,M))$
of the $NM$ pairs $(s,t)$.
\end{lem}

\begin{proof} Without loss of generality assume $N \le M$.  If for
  some pair $(s,t)$ we have $E_\f(t)$ and $E_\g(s)$ are associated to
  the same weight-$2$ cuspidal newform, then the Rankin-Selberg
  convolution will be imprimitive (and divisible by $\zeta(s)$); call
  such a pair \emph{bad}. If two elliptic curves are isomorphic, then
  they have the same $j$-invariant. Thus for a bad pair, \be j_\f(t) \
  = \ \frac{ 3A_1(t)^3}{4 A_1(t)^3 + 27 B_1(t)^2} \ = \ \frac{
    3A_2(s)^3}{4 A_2(s)^3 + 27 B_2(s)^2} \ = \ j_\g(s). \ee As we are
  assuming $j_\f(t)$ and $j_\g(s)$ are non-constant, for each fixed
  $t$ there are only finitely many solutions to $j_\g(s) = j_\f(t)$
  (the number is bounded by the degrees of $A_2(s)^3$ and $4 A_2(s)^3
  + 27 B_2(s)^2$). Thus of the $NM$ pairs $(s,t)$, at most $O(N)$ of
  the pairs have a Rankin-Selberg convolution divisible by the Riemann
  zeta function. As the only non-primitive $L$-functions $L(s,f \times
  g)$ for $f, g$ primitive weight-$2$ cuspidal newforms of levels
  $N_1$ and $N_2$ arise when $f=g$, the remaining pairs yield
  primitive $L$-functions.
\end{proof}

We now prove our main result about convolving two families of
elliptic curves.

\begin{thm}\label{thm:maintwistecf}
Consider two one-parameter families of elliptic curves (elliptic
surfaces over $\Q$): \bea \mathcal{E}_\f:& y^2
\ = \  x^3 + A_1(T)x + B_1(T) \nonumber\\
\mathcal{E}_\g:& y^2  \ = \  x^3 + A_2(S)x + B_2(S). \eea Let $\f_N$
be the specialization of $\mathcal{E}_\f$ with $t\in [N,2N-1]$,
$\g_M$ be the specialization of $\mathcal{E}_\g$ with $s \in
[M,2M-1]$, and set $\f = \cup \f_N$ and $\g = \cup \g_M$. Assume
$\log N \ll \log M \ll \log N$ and \ben
\item the first family is an elliptic curve over $\Q(T)$ of rank
$r_\f$ and non-constant $j_\f(T)$;
\item the second family is an elliptic curve over $\Q(S)$ of rank $r_\g$ and
non-constant $j_\g(S)$; \item the average log-conductor of $\f_N
\times \g_M$ satisfies \eqref{eq:lowerboundsrnmcond}; \item the
Fourier coefficients of each family satisfy either
\eqref{eq:onelevelsumbfpb} or \eqref{eq:onelevelsumbfp}. \een Then
Theorem \ref{thm:ntgRSconv} holds for the family $\f \times \g$; the
symmetry is symplectic and the rank is $0$.
\end{thm}

\begin{rek}
Rosen and Silverman \cite{RoSi} show that \eqref{eq:onelevelsumbfpb}
is a consequence of Tate's conjecture \cite{Ta}: \emph{Let
$\mathcal{E}/ \Q$ be an elliptic surface and $L_2(\mathcal{E},s)$ be
the $L$-series attached to $H^2_{\mbox{{\'e}{\rm t}}}(\mathcal{E}/
\overline{\Q}, \Q_l)$. $L_2(\mathcal{E},s)$ has a meromorphic
continuation to $\C$ and $-\mbox{ord}_{s=1} L_2(\mathcal{E},s)$ $=
\mbox{rank}\ NS(\mathcal{E}/ \Q)$, where $NS(\mathcal{E}/ \Q)$ is
the $\Q$-rational part of the N{\'e}ron-Severi group of
$\mathcal{E}$. Further, $L_2(\mathcal{E},s)$ does not vanish on the
line $\mbox{Re}(s) = 1$.}

Tate's conjecture is known for rational surfaces\footnote{An
elliptic surface $y^2 = x^3 + A(T)x + B(T)$ is rational if and only
if one of the following is true: \emph{either} $0 < \max\{3
\deg A, 2\deg B\} < 12$ \emph{or} $3\deg A =
2\deg B = 12$ and $\mbox{ord}_{T=0}T^{12} \Delta(T^{-1}) =
0$.}. Theorem \ref{thm:maintwistecf} should be true for families
with constant $j$-invariants; however, for such families Michel's
result on the average second moments of the Fourier coefficients is
not available, and one must show by direct calculation that
\eqref{eq:3iitwo} holds.
\end{rek}

\begin{proof}
For the $L$-function attached to $E_\f(t) \times E_\g(s)$, the
explicit formula \eqref{eq:14} becomes
\begin{eqnarray}\label{eq:14later}
  & & \sum_\ell \phi\left(\gamma_{E_{\f}(t)\times E_\g(s),\ell}\frac{\log R_{N,M}}{2\pi}\right) \ =
  \
  \frac{A_{E_{\f}(t)\times E_\g(s)}}{\log R_{N,M}}\hphi(0) \nonumber\\ & &
  \ \ \ \ \ \ \ \ \
  - \
  2\sum_{p}\sum_{\nu=1}^\infty\hphi\left(\frac{\nu\log p}{\log R_{N,M}}\right)
  \frac{b_{E_{\f}(t)\times E_\g(s)}(p^\nu)\log p}{p^{\nu/2}\log R_{N,M}},
\end{eqnarray} where \be A_{E_{\f}(t)\times E_\g(s)} \ = \ \log
Q(E_{\f}(t), E_\g(s)) + o(1).\ee The $o(1)$ error follows from
Theorem \ref{thm:sizearithcondprod}, where we showed $R_{N,M}$
cannot be too small. Note \eqref{eq:14later} may be slightly off in
that, if $E_\f(t) = E_\g(s)$, then the $L$-function associated to
$E_\f(t) \times E_\g(s)$ is imprimitive.
We would have a superposition of zeros of two primitive
$L$-functions, one of which has a pole. Fortunately, by Lemma
\ref{lem:nozetamosttime}, this occurs for at most $O(\min(N,M))$ of
the $NM$ pairs; as we divide by $NM$ this contribution is
negligible.\footnote{Here we are using Corollary \ref{cor:cortomainthmonconvolvingfamilies}, which says it suffices to show there is a power savings in the number of bad pairs. Alternatively, instead of using Lemma \ref{lem:nozetamosttime} we could show that the multiplicity of any elliptic curve in our parametrizations is $O(1)$.}

Thus summing \eqref{eq:14later} over $t \in [N,2N-1]$ and $s \in
[M,2M-1]$, and recalling the definition of the $1$-level density and
the average log-conductor, we find \bea\label{eq:14laterb} D_{1,\f_N
\times \g_M}(\phi) \ = \ \hphi(0) + o(1) - \frac{2}{NM}
\sum_{t=N}^{2N-1} \sum_{s=M}^{2M-1}
\sum_{p}\sum_{\nu=1}^\infty\hphi\left(\frac{\nu\log p}{\log
R_{N,M}}\right) \frac{b_{E_{\f}(t)\times E_\g(s)}(p^\nu)\log
p}{p^{\nu/2}\log R_{N,M}}.\nonumber\\ \eea

By Lemma \ref{lemftimesg}, $b_{E_{\f}(t)\times E_\g(s)}(p^\nu) =
b_{E_{\f}(t)}(p^\nu) \cdot b_{E_\g(s)}(p^\nu)$ when
$L(s,E_{\f}(t)\times E_\g(s))$ is primitive. We use this for all
$E_{\f}(t)\times E_\g(s)$, as the $O(\min(N,M))$ instances where
this is false lead to a difference that is $o(1)$.

There is trivially no contribution in \eqref{eq:14laterb} for $\nu
\ge 3$. As for each $E_{\f}(t)\times E_\g(s)$ the conductor is at
most $(NM)^b$ for some $b$, at primes dividing the conductor if
necessary we may adjust the coefficients at $p$ and $p^2$ and
introduce an error at most $o(1)$. This is because the worst case is
if $(NM)^b$ is the product of the first $\ell$ primes, where $p_\ell
\ll \log (NM)^b$. This would lead to a sum bounded by \be
\frac1{\log R_{N,M}} \sum_{p \le \log (NM)^b} \frac1{\sqrt{p}} \ \ll
\ \frac{\sqrt{\log ((NM)^b)}}{\log R_{N,M}} \ = \ o(1), \ee where
the last inequality follows from the lower bound for the average
log-conductor.

The proof is completed by showing our family is $\ntg$. We must
check the four conditions of Definition \ref{defntg}. The first is
straightforward (with Lemma \ref{lem:nozetamosttime} a key ingredient), the second (on the size of the log-conductors) follows from
our assumption that the average log-conductor satisfies
\eqref{eq:lowerboundsrnmcond}. The fourth is an easy consequence of
the Hasse bound. We are left with the third condition, which
concerns the sums over primes and squares of primes. We handle the
prime sums first.

The needed result for the sum of the Fourier coefficients at the
primes is true because \eqref{eq:onelevelsumbfp} is satisfied with
$r = 0$. To see this, note \bea\label{eq:tempecfsum1} & &
\frac{1}{NM} \sum_{t=N}^{2N-1} \sum_{s=M}^{2M-1} \sum_{p}\ \hphi\left(\frac{\log p}{\log R_{N,M}}\right)
\frac{b_{E_{\f}(t)\times E_\g(s)}(p)\log p}{\sqrt{p}\log R_{N,M}}
\nonumber\\ &  = \ & \frac{1}{NM} \sum_{t=N}^{2N-1}
\sum_{s=M}^{2M-1} \sum_{p}\ \hphi\left(\frac{\log p}{\log
R_{N,M}}\right) \frac{b_{E_{\f}(t)}(p) b_{E_{\g}(s)}(p)\log
p}{\sqrt{p}\log R_{N,M}} \nonumber\\ & = \ & \sum_{p}\ \hphi\left(\frac{\log p}{\log R_{N,M}}\right) \left[ \frac1N
\sum_{t=N}^{2N-1} b_{E_{\f}(t)}(p)\right] \left[\frac1M
\sum_{s=M}^{2M-1} b_{E_{\g}(s)}(p)\right] \frac{\log p}{\sqrt{p}\log
R_{N,M}}. \nonumber\\  \eea We analyze the $t$-sum; the $s$-sum
follows similarly. Let $a_{E_\f(t)}(p) = b_{E_\f(t)} \sqrt{p}$; by
Hasse's bound we have $|a_{E_\f(t)}(p)| \le 2 \sqrt{p}$, and these
correspond to the associated $L$-function having functional equation
$u \to 2 - u$. Let $A_p(\mathcal{E}_\f) = \frac1p \sum_{t \bmod p}
a_{E_\f(t)}(p)$. We have \bea \frac1N \sum_{t=N}^{2N-1}
b_{E_{\f}(t)}(p) & \ = \ & \frac1N\left( \frac{N}{p} \sum_{t \bmod
p} \frac{a_{E_\f(t)}(p)}{\sqrt{p}} + O(p) \right) \ = \
\frac{A_p(\mathcal{E}_\f)}{\sqrt{p}} + O\left(\frac{p}{N}\right).
 \nonumber\\
\eea The $O(p/N)$ term (and the corresponding $O(p/M)$ term from the
$s$-sum) lead to $o(1)$ contributions if $\hphi$ has suitably
restricted support. We are left with the $A_p(\mathcal{E}_\f)
A_p(\mathcal{E}_\g)/p$ term. Thus \eqref{eq:tempecfsum1} becomes \be
\sum_{p}\ \hphi\left(\frac{\log p}{\log R_{N,M}}\right)
\frac{A_p(\mathcal{E}_\f) A_p(\mathcal{E}_\g) \log p}{p^{3/2}\log
R_{N,M}} + o(1). \ee As $A_p(\mathcal{E}_\f)$ and
$A_p(\mathcal{E}_\g)$ are bounded independent of $p$ (see \cite{De},
or \cite{Mic} for an explicit bound in terms of the curves), the
above sum is $O(1)$ and hence negligible upon division by $NM$.

We are left with showing that \eqref{eq:3iione} (the second part of
the third condition of Definition \ref{defntg}) holds, i.e.,
analyzing the prime square sums (the sums of $b_{E_{\f}(t)}(p^2)$
over $t$ and $b_{E_{\g}(s)}(p^2)$ over $s$). As we have assumed
$j_\f(T)$ and $j_\g(S)$ are non-constant, this follows immediately
from work of Michel \cite{Mic}, who showed that for a one-parameter
family $\f$ over $\Q(T)$ with non-constant $j_\f(T)$ that
\begin{eqnarray}
\sum_{t\bmod p} a_{E_{\f}(t)}(p)^2\ =\ p^2 + O(p^{3/2}).
\end{eqnarray} The exponent in the error term cannot be improved in
general, and may be related to family specific lower order
correction terms to the $1$-level density; see \cite{Mil3}. From
\eqref{eq:goodeqbfp2} and our normalizations\footnote{Remember
  $b_{E_{\f}(t)}(p)\sqrt{p} = a_{E_{\f}(t)}(p)$. We must be careful in
  our normalizations, since we wish our elliptic curve $L$-functions to have a functional equation as $u \to 1-u$ (\emph{not} as $u \to
  2-u$).} we have $b_{E_{\f}(t)}(p^2) = p^{-1}a_{E_{\f}(t)}(p)^2 - 2$,
which implies \bea \sum_{t=N}^{2N-1} b_{E_{\f}(t)}(p^2) & \ = \ &
\frac{N}{p} \sum_{t \bmod p} b_{E_{\f}(t)}(p^2) + O(p) \nonumber\\ &
\ = \ & \frac{N}{p}
\sum_{t\bmod p} \frac{a_{E_{\f}(t)}(p)^2}{p} - 2N + O(p)\nonumber\\
& \ = \ & -N + O(p^2) \ = \ -|\f_N| + O(p).\eea Thus
\eqref{eq:3iione} holds with $c_\F = -1$; an analogous result holds
for sums of $b_{E_{\g}(s)}(p^2)$.

Therefore the two families have orthogonal symmetry (as was already
known), but the convolution family has symplectic symmetry
($c_{\f\times\g} = c_\f \cdot c_\g = (-1)^2 = 1$).
\end{proof}

\begin{rek}
The conditions of Theorem \ref{thm:maintwistecf} are quite weak, and
are easily seen to be satisfied in many cases of interest (for
example, by many of the families studied in \cite{ALM,Fe}).
\end{rek}

\appendix


\bigskip
\section{Average Log-Conductors for Elliptic
Curve Families}\label{sec:proofthmsizearithcondprod}

We prove Theorem \ref{thm:sizearithcondprod}. The upper bound
follows trivially from \eqref{eq:sizecondsc1c2} and bounds relating
the discriminant of an elliptic curve to its conductor. We prove the
lower bound through a series of lemmas. We first introduce some
notation. Let $f_1(x), \dots, f_{k_1}(x)$ be the distinct monic
irreducible factors of $\Delta_\f(x) c_{\f,4}(x) c_{\f,6}(x)$, and
let $g_1(x), \dots, g_{k_2}(x)$ be the distinct monic irreducible
factors of $\Delta_\g(x) c_{\g,4}(x) c_{\g,6}$. By relabeling if
necessary, we may assume $f_1(x)|\Delta_\f(x)$, and similarly
$g_1(x)|\Delta_\g(x)$.

Further, we may assume \emph{all} of the $k_1+k_2$ polynomials are
relatively prime. If some $f_i(x)$ and $g_j(x)$ were not relatively
prime, then we could find a fixed $x'$ such that $g_1(x+x'), \dots,
g_{k_2}(x+x')$ are relatively prime (as functions of $x$) to the
$f_i(x)$'s. Thus instead of considering the interval $[M,2M]$ we
would consider the interval $[M-c, 2M - c]$, which for $M$ large is
still approximately $[M,2M]$. Hence we may assume $f_1(x) \notdiv
\Delta_\g(x)$ and $g_1(x) \notdiv \Delta_\f(x)$, and without loss of
generality we may assume $N \le M$.

The proof is completed by showing that, for some constants $\gep,
\delta, a > 0$, at least $\gep NM / (\log \log N)^a$ of the
$L(u,E_\f(t) \times E_\g(s))$ have arithmetic conductor at least
$(NM)^\delta$. We do this by showing for at least $\gep NM / (\log
\log N)^a$ of the pairs $(t,s)$ that we can find a number at least
$N^\delta$ dividing $f_1(t)$ and the conductor $C_\f(t)$ but not any
other $f_i(t)$, and a number at least $M^\delta$ (relatively prime
to the number at least $N^\delta$) dividing $g_1(s)$ and $C_\g(s)$
but not any other $g_i(s)$. As the arithmetic conductor is an
integer, it will then be divisible by at least $(NM)^\delta$, which
implies the lower bound in \eqref{eq:lowerboundsrnmcond}. We do this
through the following series of lemmas.

\begin{lem}\label{lem:condsize1}
There exists an integer $c$ (a product of distinct primes) and an
integer $r$ such that, for $i \neq j$, for a positive fraction of $t
\in [N,2N]$ we have $(f_i(t),f_j(t)) | c^r$ (and similarly for the
$g$'s).
\end{lem}

\begin{proof} Let $i \neq j$. By the Euclidean algorithm, there exists a
$c_{ij}$ (independent of $x$) such that if $p|(f_i(x),f_j(x))$ then
$p|c_{ij}$. Let $c_f$ be the product of $6$ and the prime divisors
$p \ge 5$ of the $c_{ij}$'s. Choose an $x_0$ such that $f_i(x_0)
\neq 0$ for all $i$. Let $r_f$ be the largest integer such that if
$p|c_f$ then $p^{r+1}$ divides none of the $f_i(x_0)$. Then for all
$i \neq j$ the greatest common divisor of $f_i(c_f^{r_f+1}x + x_0)$
and $f_j(c_f^{r_f+1}x + x_0)$ divides $c_f^{r_f}$. We similarly
construct $c_g$ and $r_g$ so that the greatest common divisors of
the $g_i$'s divides $c_g^{r_g}$. Let $c$ equal the product of the
prime divisors of $c_f c_g$ and $r = \max (r_f,r_g)$. We change
variables, sending $t \to c^{r+1} t + x_0$. A positive fraction of
$t \in [N,2N]$ satisfy this condition. We similarly change $s \to
c^{r+1} s + s_0$. For ease of exposition we denote these polynomials
by $\tf_i$ and $\tg_j$; thus $\tf_i(x) = f_i(c^{r+1}x + x_0)$.
\end{proof}

It is possible that $\tf_i(x)$ is divisible by a fixed square (or
higher power) for all $x$; for example, $x^4 - x^2 + 20$ is always
divisible by $4$. Further, if $\Delta_\f(x)c_{\f,4}(x)c_{\f,6}(x) =
a_\f \prod_i f_i(x)$, for some $x$ an $\tf_i(x)$ could share a
factor with $a_\f$. The following lemma handles such primes.

\begin{lem}\label{lem:condsize2} Notation as in Lemma
\ref{lem:condsize1}, let $C$ be the product of all numbers that
divide an $\tf_i(x)$ for all $x$, a $\tg_j(x)$ for all $x$, or $a_\f
a_\g$. Then there exists an integer $m$ such that, for a positive
fraction of $t \in [N,2N]$, Lemma \ref{lem:condsize1} holds
\emph{and} if $p|C$ then $p^m$ does not divide $f_i(t)$ for all $i$
(and similarly for the $g$'s).
\end{lem}

\begin{proof} Let $\tf_i$ and $\tg_j$ be as in Lemma \ref{lem:condsize1}.
Choose an $x_1$ such that $\tf_i(x_1)$ and $\tg_j(x_1)$ are non-zero
for all $i$ and $j$. Arguing as before, after a simple linear change
of variables we can ensure that at most a fixed power of $C$ divides
our polynomials for any $x$. Specifically, consider $\tf_i(C^m x +
x_1)$; for $m$ sufficiently large, if $p|C$ then $p^m \notdiv
\tf_i(C^m x + x_1)$ (and the same is true for the $\tg_j$'s). Let
$\whf_i(x) = \tf_i(C^m x + x_1)$ (and similarly for $\whg_j$).
\end{proof}

We have shown that for a positive fraction of all $t \in [N,2N]$ and
$s \in [M,2M]$: (i) the greatest common divisor of the $\whf_i(t)$'s
is at most $c^r$ and the greatest common divisor of the
$\whg_j(s)$'s is at most $c^r$; (ii) the product of all the squares
or factors of $a_\f a_\g$ that divide a $\whf_i(t)$ for all $t$ is
at most $C^m$ (and similarly for $\whg_j(s)$). We would like to say
the arithmetic conductor is at least $\whf_1(t) \whg_1(s)$, except
there are two problems: (i) we must show $\whf_1(t)$ divides the
conductor of $E_\f(t)$ (and similarly for $\whg_1(s)$ and
$E_\g(s)$); (ii) we must show $(\whf_1(t), \whg_1(s))$ is small. We
handle (i) first.

\begin{lem}\label{lem:condsize3} Let $d= \max(\sum_i \deg f_i, \sum_j \deg g_j)
+2$. For a positive fraction of $t \in [N,2N]$ and $s\in [M,2M]$,
the results of Lemmas \ref{lem:condsize1} and \ref{lem:condsize2}
hold, the conductor of $E_\f(t)$ is $\gg N^{1/d}$ and the conductor
of $E_\g(s)$ is $\gg M^{1/d}$
\end{lem}

\begin{proof} Notation as in Lemmas \ref{lem:condsize1} and
\ref{lem:condsize2} and conditions as in Theorem
\ref{thm:sizearithcondprod}, we show that for a positive fraction of
the time that the conductor of $E_\f(t)$ is $\gg N^{1/d}$. Recall
the following basic facts (see for example \cite{Nag}) for an
integral polynomial $D(t)$ of degree $k$ and discriminant $\delta$:

\begin{enumerate}
  \item Let $p$ be a prime not dividing the coefficient of
$x^k$. Then $D(t) \equiv 0$ mod $p$ has at most $k$ incongruent
solutions.
  \item Suppose $p \notdiv \delta$. Then the number of incongruent solutions
  of $D(t) \equiv 0$ mod $p$ equals the number of incongruent solutions of
$D(t) \equiv 0$ mod $p^{\alpha}$.
\end{enumerate}

Note that if the discriminant of $h(x)$ is $\delta$, then the
discriminant of $h(ax+b)$ is $a^n \delta$ for some $n$. Let $D$ be
the product of the prime divisors of the discriminants and leading
coefficients of all the $\whf_i$'s and $\whg_j$'s, as well as any
missing primes at most $d$. We make one last change of variables:
for sufficiently large $n$ consider \be F_i(x) \ = \ \whf_i(D^n x +
x_2), \ \ \ G_j(x) \ = \ \whg_j(D^n x + x_2), \ee where $x_2$ is
chosen so that all $\whf_i(x_2)$ and $\whg_j(x_2)$ are non-zero. The
advantage is that the degree of divisibility of $F_i(x)$ (resp.,
$G_j(x)$) by primes dividing the discriminants, $c_{\f,4}(x)$ and
$c_{\f,6}(x)$ (resp., $c_{\g,4}(x)$ and $c_{\g,6}(x)$), leading
coefficients or at most $d$ is bounded independent of $x$, say by
$k$. It is now immediate that, for a positive fraction of $x$,
$F_i(x)$ has a $d$\textsuperscript{th} power free factor $\gg
N^{\deg F_i / d}$. To see this, let $\nu_{F_i}(p^d)$ denote the
number of solutions to $F_i(x) \equiv 0 \bmod p^{d}$. For $p \notdiv
D$, $p$ does not divide the discriminant of $F_i$ and thus
$\nu_{F_i}(p^d) = \nu_{F_i}(p) \le \deg F_i$. Thus the fraction of
$t$ giving $F_i(t)$ $d$-power free (except for divisors of $D$) is
at least \be\label{eq:numdivpnotdivDFi} \prod_{p \notdiv D} \left(1
- \frac{\nu_{F_i}(p^d)}{p^d} \right) \ \ge \ \prod_{p \notdiv D}
\left(1 - \frac{\deg f_i}{p^d} \right) \ \ge \ \prod_{p \notdiv D}
\left(1 - \frac{1}{p^{d-1}} \right). \ee As $d$ was chosen to be at
least $3$, this last factor is larger than $\prod_p (1 - p^{-2}) =
6/\pi^2$. By our linear change of variables (how we defined the
$F_i$), the number of times a $p|D$ divides $F_i(x)$ is bounded
independent of $x$ and $i$, say by $k$. Thus, for a positive
fraction of $t$, $F_i(t)$ has a $d$ power free part at least
$F_i(t)^{1/d} / D^k$. As the greatest common divisors of any two of
the $F_i$ is at most $c^r$, for a positive fraction of $t$ we have
$F_1(t)$ has a $d$ power free factor of size at least $F_1(t)^{1/d}
/ c^r D^k$ that is relative prime to the $F_i(t)$ for $i \neq 1$.

We need only show this factor (which is at least $F_1(t)^{1/d} / c^r
D^k$) divides the conductor. This follows by showing the conductor
of the elliptic curve $y^2 = x^3 + A_1(t)x + B_1(t)$ is minimal for
each $p \notdiv cCD$ that divides $F_1(t)^{1/d}$. This follows from
our assumption that $j_\f(T)$ is not constant, as this implies that
$c_{\f,4}(x)$ and $c_{\f,6}(x)$ are not identically zero. Thus
neither are $C_{\f,4}(x)$ or $C_{\f,6}(x)$ (where we have used the
obvious notation to represent the linear change of variables). By
assumption (see the conditions of Theorem
\ref{thm:sizearithcondprod}), as the irreducible polynomial factors
of $c_{\f,4}(x)$ and $c_{\f,6}(x)$ were included in our list of the
$f_i$'s, and we assumed $f_1(x)$ is relatively prime to either
$c_{\f,4}(x)$ or $c_{\f,6}(x)$, for $p \notdiv cCD$ with $p|F_1(t)$,
$p$ cannot divide both $C_{\f,4}(x)$ and $C_{\f,6}(x)$. Thus the
elliptic curve is minimal for such primes $p$, implying the
conductor is at least $F_1(t)^{1/d} / c^r C^m D^k$.

Thus as $N \to \infty$, for a positive fraction of $t$ the conductor
of $E_\f(t)$ is $\gg N^{\deg F_1 / d} \gg N^{1/d}$; an analogous
statement holds a positive fraction of the time for the conductor of
$E_\g(s)$. We call such $t$ and $s$ \emph{good}.
\end{proof}

\begin{rek} We chose $d = \max(\sum_i \deg f_i, \sum_j \deg g_j) +2$
and not $\max(\max_i \deg f_i$, $\max_j \deg g_j)$ $+2$ because of
Lemma \ref{lem:finallemproofthmcondsize}.
\end{rek}

The following lemma completes the proof of Theorem
\ref{thm:sizearithcondprod}.

\begin{lem}\label{lem:finallemproofthmcondsize}
Notation as in Theorem \ref{thm:sizearithcondprod}, for some $\gep,
a > 0$ for at least $\gep NM / (\log \log N)^a$ of the pairs $(t,s)
\in [N,2N] \times [M,2M]$ the results of Lemmas \ref{lem:condsize1}
through \ref{lem:condsize3} hold, and the greatest common divisor of
$\prod_i F_i(t)$ and $\prod_j G_j(s)$ is bounded independent of $t$
and $s$.
\end{lem}

\begin{proof} Consider the positive fraction of $t$ and $s$ that are good. We
must make sure that each such $G_j(s)$ is essentially relatively
prime to the $F_i(t)$. If so, then since the arithmetic conductor is
an integer it would have to be $\gg N^{1/d} M^{1/d}$ (remember the
arithmetic conductor comes from the arithmetic conductors of
$E_\f(t)$ and $E_\g(s)$, and these are $\gg N^{1/d}$ and $\gg
M^{1/d}$). For a good $t$, the worst case for common factors of
$\prod_i F_i(t)$ and $\prod_j G_j(s)$ is when $\prod_i F_i(t)$ is
the product of the first $\ell$ primes. We can easily handle the
bounded contributions from $c$ (Lemma \ref{lem:condsize1}), $C$
(Lemma \ref{lem:condsize2}) or $D$ (see the proof of Lemma
\ref{lem:condsize3}), and thus we need only investigate primes $p$
such that $p > D$ and $p\notdiv cC$. Letting $\mu = \deg
\Delta_\f(x) C_{\f,4}(x) C_{\f,6}(x)$, the product of the $F_i(t)$'s
is $\ll N^\mu$.
Thus \bea \prod_{p \le p_\ell} p & \ \ll \ & N^\mu \nonumber\\
\sum_{p \le p_\ell} \log p & \ll & \mu \log N \ \ \Rightarrow \ \
p_\ell \ \ll \ \mu \log N. \eea

We may need to discard some good $s$ because a $G_j(s)$ is not
relatively prime to $p_1 \dots p_\ell$. Thus, even though we are
going to use $F_1(t)$ and $G_1(s)$, we must make sure that there are
no large common factors of $\prod_i F_i(t)$ and $\prod_j G_j(s)$, as
otherwise the conductor of the Rankin-Selberg convolution could be
reduced (see the division by the greatest common divisor in the left
hand side of \eqref{eq:sizecondsc1c2}). As $d > \sum_j \deg g_j$, we
have $\sum_j \mu_{G_j}(p) < d$ for $p \notdiv D$; this allows us to
obtain the needed estimate.

We have already handled $p \le d$ and $p$ dividing a discriminant in
our construction of the good $s$. For each $j$, the product in
\eqref{eq:numdivpnotdivDFi} for $G_1(s)$ is modified by a factor no
worse than \bea \prod_{d < p \le p_\ell \atop p \notdiv cCD} \left(1
- \frac{\sum_j\nu_{G_j}(p)}{p}\right) & \ \ge
\ & \prod_{d < p \le p_\ell} \left(1 - \frac{d}{p}\right) \nonumber\\
& \gg & \exp\left(-2009d\cdot\log \sum_{d < p \le p_\ell}
\frac{1}{p}\right) \nonumber\\ & \gg & \left(\log
p_\ell\right)^{-2009d} \ \gg \ \left(\log \log N\right)^{-2009d}
\eea (the last bound is Mertens' theorem, see for example
\cite{Da}).
We therefore obtain a $d$ power free factor of the conductor of
$E_\g(s)$ whose common factors with $F_1(t)\cdots F_{k_1}(t)$ is
bounded by $c^r D^k$.
\end{proof}

This completes the proof of Theorem \ref{thm:sizearithcondprod}.
\hfill $\Box$



\ \\


\begin{thebibliography}{PTW021} 

\bibitem[ALM]{ALM}
\newblock S. Arms, A. Lozano-Robledo and S. J. Miller, \emph{Constructing
one-parameter families of elliptic curves over $\Q(T)$ with moderate
rank}, Journal of Number Theory \textbf{123} (2007), no.~2, 388--402.

\bibitem[Bor]{Bor} \newblock A. Borel, \emph{Automorphic $L$-functions,}
  in \emph{Automorphic Forms, Representations, and $L$-functions.}
  Proc.~Symp.~Pure Math. \textbf{33} part~2, 21--61, Amer.\ Math.\ Soc.,
  Providence 1979.

\bibitem[BCDT]{BCDT}
\newblock C. Breuil, B. Conrad, F. Diamond and R. Taylor, \emph{On
the modularity of elliptic curves over $\Q$: wild $3$-adic
exercises}, J. Amer. Math. Soc. \textbf{14} (2001), no. 4, 843--939.

\bibitem[Bu]{Bu}
\newblock D. Bump, \emph{The Rankin-Selberg method: a survey}.
In \emph{Number theory, trace formulas and discrete groups (Oslo,
1987)}, 49--109, Academic Press, Boston, MA, 1989.

\bibitem[BH]{BH}
\newblock C. J. Bushnell and G. Henniart, \emph{On
certain dyadic representations}, appendix to [KiSh1].

\bibitem[CM]{CM}
\newblock J. Cogdell and P. Michel, \emph{On the complex moments of
  symmetric power $L$-functions at $s=1$,} Int.\ Math.\ Res.\ Notices
2004, no.\ 31, 1561-1617.

\bibitem[CFKRS]{CFKRS}
\newblock B. Conrey, D. Farmer, P. Keating, M. Rubinstein and N.
Snaith, \emph{Integral moments of $L$-functions}, Proc. London Math.
Soc. (3)  \textbf{91} (2005),  no. 1, 33--104.

\bibitem[Da]{Da}
\newblock H. Davenport, \emph{Multiplicative Number Theory, $2$nd edition},
 Graduate Texts in Mathematics \textbf{74}, Springer-Verlag, New York,
 $1980$, revised by H. Montgomery.

\bibitem[De]{De}
\newblock P. Deligne, \emph{La conjecture de Weil II}, Publ. Math.
IHES \textbf{52} (1981), 313--428.

\bibitem[DM]{DM}
\newblock E. Due\~nez and S. J. Miller, \emph{The low lying zeros of a
$\text{GL}(4)$ and a $\text{GL}(6)$ family of
$L$-functions}, Compositio Mathematica \textbf{142} (2006), no. 6, 1403--1425.

\bibitem[Fe]{Fe}
\newblock S. Fermigier, \emph{\'Etude exp\'erimentale du rang de
familles de courbes elliptiques sur $\Q$}, Exper. Math. \textbf{5}
(1996), 119--130.

\bibitem[FI]{FI}
\newblock E. Fouvry and H. Iwaniec, \emph{Low-lying zeros of dihedral
$L$-functions}, Duke Math. J. \textbf{116} (2003), no. 2, 189--217.

\bibitem[Ga]{Ga}
\newblock P. B. Garrett, \emph{Decomposition of Eisenstein series:
Rankin triple products},  Ann. of Math. (2)  \textbf{125} (1987),
no. 2, 209--235.

\bibitem[Gao]{Gao}
\newblock P. Gao, \emph{$N$-level density of the low-lying zeros of
quadratic Dirichlet $L$-functions}, Ph.~D thesis, University of
Michigan, 2005.

\bibitem[Gel]{Gel}
\newblock S.~Gelbart, \emph{Automorphic Forms on Adele Groups,} Annals of
Mathematics Studies \textbf{83}, Princeton University Press, Princeton, NJ (1975).

\bibitem[GJ]{GJ}
\newblock S. Gelbart and H. Jacquet, \emph{A relation between
automorphic representations of GL(2) and GL(3)}, Ann.~Sci.~\'Ecole
Normale Sup., 4e s\'erie, Vol.~11 (1978), 471--552.

\bibitem[Go]{Go}
\newblock D. Goldfeld, \emph{Conjectures on elliptic curves over
quadratic fields}, Number Theory (Proc. Conf. in Carbondale,
$1979$), Lecture Notes in Math. \textbf{751}, Springer-Verlag,
$1979$, 108--118.

\bibitem[GR]{GR}
\newblock I. Gradshteyn and I. Ryzhik, \emph{Tables of Integrals,
Series, and Products}, New York, Academic Press, 1965.

\bibitem[G\"u]{Gu}
\newblock A. G\"ulo\u{g}lu, \emph{Low Lying Zeros of Symmetric
Power $L$-Functions},  Int. Math. Res. Not. (2005),  no. 9,
517--550.

\bibitem[HaMi]{HaMi} \newblock G. Harcos and P. Michel, \emph{The
    subconvexity problem for Rankin-Selberg $L$-functions},
  Invent. Math.  \textbf{163} (2006), no. 3, 581--655.

\bibitem[HW]{HW}
\newblock G. Hardy and E. Wright, \emph{An Introduction to the
Theory of Numbers}, fifth edition, Oxford Science Publications,
Clarendon Press, Oxford, 1995.

\bibitem[Hej]{Hej}
\newblock D. Hejhal, \emph{On the triple correlation of zeros of
the zeta function}, Internat. Math. Res. Notices 1994, no. 7,
294--302.


\bibitem[HM]{HM}
\newblock C. Hughes and S. J. Miller, \emph{Low-lying zeros of $L$-functions
with orthogonal symmetry}, Duke Math.~J.
\textbf{136} (2007), no.~1, 115--172.

\bibitem[HR]{HR}
\newblock C. Hughes and Z. Rudnick, \emph{Linear statistics of
low-lying zeros of $L$-functions}, Quart. J. Math. Oxford
\textbf{54} (2003), 309--333.

\bibitem[Iw]{Iw}
\newblock H. Iwaniec, \emph{Introduction to the Spectral Theory of
Automorphic Forms}, 2nd edition, Graduate Studies in Mathematics
\textbf{53}, AMS, 2002.

\bibitem[ILS]{ILS}
\newblock H. Iwaniec, W. Luo and P. Sarnak, \emph{Low lying zeros of
families of $L$-functions}, Inst. Hautes \'Etudes Sci. Publ. Math.
\textbf{91} (2000), 55--131.

\bibitem[Jac]{Jac}
\newblock H. Jacquet, \emph{Principal $L$-functions of the
    Linear Group,} in \emph{Automorphic Forms, Representations, and
    $L$-functions.}  Proc.~Symp.~Pure Math. \textbf{33} part~2, 63--86,
  Amer.\ Math.\ Soc., Providence 1979.

\bibitem[JPS]{JPS} \newblock H. Jacquet, I. I. Piatetski-Shapiro and J.
  A. Shalika, \emph{Rankin Selberg convolutions.} Amer.\ Jour.\ of Math.
  \textbf{105} (1983), 367--464.

\bibitem[JS]{JS}
\newblock H. Jacquet and J. A. Shalika, \emph{On Euler products and the classification of automorphic representations I}, Amer. Jour. of Math. \textbf{103} (1981), 499--558.

\bibitem[KaSa1]{KaSa1}
\newblock N. Katz and P. Sarnak, \emph{Random Matrices, Frobenius
Eigenvalues and Monodromy}, AMS Colloquium Publications \textbf{45},
AMS, Providence, 1999.

\bibitem[KaSa2]{KaSa2}
\newblock N. Katz and P. Sarnak, \emph{Zeros of zeta functions and symmetries},
Bull. AMS \textbf{36} (1999), 1--26.


\bibitem[KeSn]{KeSn}
\newblock J. P. Keating and N. C. Snaith, \emph{Random
matrices and $L$-functions}, Random matrix theory, J. Phys. A
\textbf{36} (2003), no. 12, 2859--2881.

\bibitem[K]{K} \newblock H. Kim, \emph{Functoriality for the exterior
square of ${\rm GL}\sb 2$ and the symmetric fourth of ${\rm GL}\sb
2$} (with appendix 1 by D. Ramakrishnan and appendix 2 by H. Kim and P. Sarnak), Jour. AMS \textbf{16} (2003), no. 1, 139--183.

\bibitem[KiSh1]{KiSh1}
\newblock H. Kim and F. Shahidi, \emph{Functorial products for
${\rm GL}\sb 2\times{\rm GL}\sb 3$ and the symmetric cube for ${\rm
GL}\sb 2$}, Annals of Math. \textbf{155} (2002), 837--893.

\bibitem[KiSh2]{KiSh2}
\newblock H. Kim and F. Shahidi, \emph{Cuspidality of symmetric powers
  with applications},  Duke Math. J. \textbf{112} (2002), no.~1, 177-197.

\bibitem[Kn]{Kn}
\newblock A. W. Knapp, \emph{Local Langlands Correspondence: the
  archimedean case,} in \emph{Motives,} Proc.\ Symp.\ Pure
Math. \textbf{55}, part 2 (1994).

\bibitem[LS]{LS}
\newblock W. Luo and P. Sarnak, \emph{Mass equidistribution for Hecke eigenforms},
Comm. Pure Appl. Math. \textbf{56} (2003), no. 7, 874--891.


\bibitem[Mic]{Mic}
\newblock P. Michel, \emph{Rang moyen de familles de courbes elliptiques
et lois de Sato-Tate}, Monat. Math. \textbf{120} (1995), 127--136.

\bibitem[Mil1]{Mil1}
\newblock S. J. Miller, \emph{$1$- and $2$-Level Densities for Families of Elliptic
Curves: Evidence for the Underlying Group Symmetries}, P.H.D.
Thesis, Princeton University, 2002.
\texttt{http://www.williams.edu/go/math/sjmiller/public\underline{\ \ }html/}

\bibitem[Mil2]{Mil2}
\newblock S. J. Miller, \emph{$1$- and $2$-level densities for families of elliptic
curves: evidence for the underlying group symmetries}, Compositio
Mathematica \textbf{140} (2004), 952--992.

\bibitem[Mil3]{Mil3}
\newblock S. J. Miller, \emph{Variation in the number of points on elliptic curves and
applications to excess rank}, C. R. Math. Rep. Acad. Sci. Canada
\textbf{27} (2005), no. 4, 111--120.

\bibitem[Mil4]{Mil4}
S. J. Miller (with an appendix by E. Due\~nez), \emph{Investigations of zeros near the central point
of elliptic curve $L$-functions}, Experimental Mathematics \textbf{15} (2006), no. 3, 257--279.

\bibitem[Mil5]{Mil5}
S. J. Miller, \emph{A symplectic test of the $L$-Functions Ratios Conjecture}, Int Math Res Notices (2008) Vol. 2008, article ID rnm146, 36 pages, doi:10.1093/imrn/rnm146.

\bibitem[Mil6]{Mil6}
\newblock S. J. Miller, \emph{Lower
order terms in the $1$-level density for families of holomorphic
cuspidal newforms},  Acta Arithmetica \textbf{137} (2009), 51--98.

\bibitem[Mil7]{Mil7}
\newblock S. J. Miller, \emph{An
orthogonal test of the $L$-functions Ratios Conjecture}, Proceedings of the London Mathematical Society 2009, doi:10.1112/plms/pdp009.

\bibitem[Mil8]{Mil8}
\newblock S. J. Miller, \emph{Density functions for families of
Dirichlet characters}, preprint.


\bibitem[Mon]{Mon}
\newblock H. Montgomery, \emph{The pair correlation of zeros of the zeta
function}, Analytic Number Theory, Proc. Sympos. Pure Math.
\textbf{24}, Amer. Math. Soc., Providence, 1973, 181--193.

\bibitem[Nag]{Nag}
\newblock T. Nagell, \emph{Introduction to Number Theory}, Chelsea
Publishing Company, New York, $1981$.

\bibitem[Od1]{Od1}
\newblock A. Odlyzko, \emph{On the distribution of spacings
between zeros of the zeta function}, Math. Comp. \textbf{48} (1987),
no. 177, 273--308.

\bibitem[Od2]{Od2}
\newblock A. Odlyzko, \emph{The $10^{22}$-nd zero of the Riemann zeta function}, Proc.
Conference on Dynamical, Spectral and Arithmetic Zeta-Functions, M.
van Frankenhuysen and M. L. Lapidus, eds., Amer. Math. Soc.,
Contemporary Math. series, 2001,
http://www.research.att.com/$\sim$amo/doc/zeta.html

\bibitem[Ram]{Ram} \newblock D. Ramakrishnan, \emph{Modularity of the
Rankin-Selberg $L$-series, and multiplicity one for ${\rm SL}(2)$},
Ann.~of Math. (2) \textbf{152} (2000), no. 1, 45--111.

\bibitem[RR]{RR}
\newblock G. Ricotta and E. Royer, \emph{Statistics for low-lying
zeros of symmetric power $L$-functions in the level aspect},
preprint. \texttt{http://arxiv.org/abs/math/0703760}


\bibitem[RoSi]{RoSi}
\newblock M. Rosen and J. Silverman, \emph{On the rank of an elliptic
surface}, Invent. Math. \textbf{133} (1998), 43--67.

\bibitem[Ro]{Ro}
\newblock E. Royer, \emph{Petits z\'{e}ros de fonctions $L$
de formes modulaires}, Acta Arith. \textbf{99} (2001), 47--172.

\bibitem[Rub]{Rub}
\newblock M. Rubinstein, \emph{Low-lying zeros of
$L$-functions and random matrix theory}, Duke Math. J. \textbf{109}
(2001), no. 1, 147--181.

\bibitem[RS]{RS}
\newblock Z. Rudnick and P. Sarnak, \emph{Zeros of principal $L$-functions
 and random matrix theory},  Duke Math. J. \textbf{81}
(1996), 269--322.

\bibitem[St]{St}
J. Stopple, \emph{The quadratic character experiment}, preprint. \texttt{http://arxiv.org/abs/0802.4255}

\bibitem[Ta]{Ta}
\newblock J. Tate, \emph{Algebraic cycles and the pole of zeta
functions}, Arithmetical Algebraic Geometry, Harper and Row, New
York, $1965$, $93-110$.

\bibitem[TW]{TW}
R. Taylor and A. Wiles, \emph{Ring-theoretic properties of certain
Hecke algebras}, Ann. Math. \textbf{141} (1995), 553--572.

\bibitem[Wi]{Wi}
A. Wiles, \emph{Modular elliptic curves and Fermat's last theorem},
Ann. Math. \textbf{141} (1995), 443--551.

\bibitem[Yo1]{Yo1}
\newblock M. Young, \emph{Lower-order terms of the 1-level density of families of elliptic
curves},  Int. Math. Res. Not. (2005),  no. 10, 587--633.

\bibitem[Yo2]{Yo2}
\newblock M. Young, \emph{Low-lying zeros of families of elliptic curves},
J. Amer. Math. Soc. \textbf{19} (2006), no. 1, 205--250.

\end{thebibliography}
\end{document}